\documentclass[11pt]{article}
\usepackage{amsmath,amsthm,amssymb}
\usepackage[dvips]{graphicx}
\usepackage[mathscr]{eucal}
\usepackage[margin=1in]{geometry}
\usepackage{standardmath}
\usepackage{psfrag}

\newcommand{\olist}[1]{({#1})^{\circ}}
\newcommand{\amap}[1][M]{{\mathscr{#1}}}
\newcommand{\aset}[1]{\mathscr{#1}}
\newcommand{\acore}[1][M]{\mathscr{#1}^c}

\newcommand{\IPhur}[2]{\smash{\widetilde{H}_{#1}(\vec{#2})}}
\newcommand{\IPshur}[2]{\smash{\widetilde{M}_{#1}(\vec{#2})}}
\newcommand{\IPshurGS}[1]{\smash{\widetilde{\Psi}_{#1}}}
\newcommand{\IPsmhur}[3]{\smash{\widetilde{S}_{#1}(\vec{#2}\,;\,\vec{#3})}}
\newcommand{\IPsmhurGS}[1]{\smash{\widetilde{\Gamma}_{#1}}}

\newcommand{\IPkshurGS}[2]{\smash{\widetilde{\Psi}_{#1,#2}}}

\newcommand{\Phur}[2]{H_{#1}(\vec{#2})}
\newcommand{\Pshur}[2]{M_{#1}(\vec{#2})}
\newcommand{\PshurGS}[1]{\Psi_{#1}}
\newcommand{\Psmhur}[2]{S_{#1}(\vec{#2})}
\newcommand{\PsmhurGS}[1]{\Gamma_{#1}}
\newcommand{\Pkhur}[2]{H_{#1}(#2)}
\newcommand{\PkshurGS}[2]{\Psi_{#1,#2}}

\renewcommand{\vec}[1]{\mathbf{#1}}
\newcommand{\rt}{w}
\newcommand{\irt}{\tilde{w}}
\newcommand{\Dset}[2]{\mathbb{N}_{#1}^{#2}}
\renewcommand{\bold}[1]{\textbf{#1}}

\title{Minimal Transitive Factorizations of Permutations into Cycles}
\author{John Irving}
\date{\today}

\begin{document}

\maketitle
\begin{abstract}
We introduce a new approach to an enumerative problem closely linked with the geometry of branched coverings;
that is, we study the number $\Pkhur{\a}{i_2,i_3,\ldots}$ of ways a given permutation (with cycles described
by the partition $\a$) can be decomposed into a product of exactly $i_2$ 2-cycles, $i_3$ 3-cycles,
\emph{etc.}, with certain minimality and transitivity conditions imposed on the factors. The method is to
encode such factorizations as planar maps with certain \emph{descent structure} and apply a new combinatorial
decomposition to make their enumeration more manageable. We apply our technique to determine
$\Pkhur{\a}{i_2,i_3,\ldots}$ when $\a$ has one or two parts, extending earlier work of Goulden and Jackson.
We also show how these methods are readily modified to count \emph{inequivalent} factorizations, where
equivalence is defined by permitting commutations of adjacent disjoint factors.  Our technique permits a
substantial generalization of recent work of Goulden, Jackson, and Latour, while allowing for a considerable
simplification of their analysis.
\end{abstract}

\section{Introduction}

We begin with a brief review of standard notation. We write $\a \ptn n$ (respectively, $\a \cmp n$) to
indicate that $\a$ is a partition (composition) of $n$, and denote by $\len{\a}$ the number of parts of $\a$.
The partition with $m_i$ parts equal to $i$ is denoted by $[1^{m_1} 2^{m_2} \cdots]$.

The \emph{cycle type} of a permutation $\p$ in the symmetric group $\Sym{n}$ is the partition of $n$
determined by the lengths of the disjoint cycles comprising $\p$. The conjugacy class of $\Sym{n}$ consisting
of all permutations of cycle type $\l \ptn n$ is denoted by $\Class{\l}$.  This notation is extended to allow
$\l$ to be a composition, in which case the ordering of the parts of $\l$ is simply ignored.  Members of
$\Class{[1^{n-k}\,k]}$ are called \emph{cycles of length $k$}, or \emph{$k$-cycles}, while elements of
$\Class{[n]}$ are \emph{full cycles} in $\Sym{n}$.

For vectors $\mathbf{j}=(j_1,\ldots,j_m)$ and
$\mathbf{x}=(x_1,\ldots,x_m)$ we use the abbreviations
$\mathbf{x}^\mathbf{j} = x_1^{j_1} \cdots x_m^{j_m}$ and
$\mathbf{j}! = j_1! \cdots j_m!$. Finally, if $f \in
\Q[[\mathbf{x}]]$ is a formal power series, then we write
$[\mathbf{x}^{\mathbf{j}}]\,f(\mathbf{x})$ for the coefficient of
the monomial $\mathbf{x}^{\mathbf{j}}$ in $f(\mathbf{x})$.

\subsection{Factorizations of Permutations into Cycles}
\label{ssec:factorizations}

A \bold{transitive factorization} of a permutation $\p \in \Sym{n}$
is a tuple $F=(\s_r,\ldots,\s_1)$ of permutations $\s_i \in \Sym{n}$
such that (1) $\p = \s_r \cdots \s_1$, and (2) the group
$\ang{\s_1,\ldots,\s_r}$ generated by the factors acts transitively
on $\Sym{n}$. If $\s_i \in \Class{\b_i}$ for $1 \leq i \leq r$ and
$\p \in \Class{\a}$, then one can show~\cite{gj-transitive} that
\begin{equation}
\label{eq:genusdefn}
    nr - \sum_{i=1}^r \len{\b_i} \geq  n + \len{\a} - 2.
\end{equation}
In the case of equality above, $F$ is said to be \bold{minimal transitive}. For example, since
\begin{align}
    \label{eq:cycexmp}
    (1\,2\,3\,4\,5)(6\,7\,8)(9)
    &=
    (3\,4\,7\,9) \cdot
    (6\,9\,7\,8) \cdot
    (1\,2) \cdot
    (2\,4\,5) \cdot (2\,6),
\end{align}
the tuple $((3479),(6978),(12),(245),(26))$ is a factorization of
$(12345)(678)(9)$. (Fixed points have been suppressed in the
factors.) This factorization is easily verified to be minimal
transitive.

Minimal transitive factorizations have been very well studied, with much of this attention stemming from the
fact that they serve geometers as combinatorial models for branched coverings of the sphere by the sphere. In
this context, transitivity guarantees connectedness of the associated covering, while minimality implies the
covering surface is the sphere. For further information on these connections see~\cite{melou-schaeffer} and
the references therein. 

The focus of this paper is the class of \bold{cycle factorizations},
by which we mean minimal transitive factorizations, such
as~\eqref{eq:cycexmp} above, whose factors are all cycles of length
at least two.\footnote{This condition avoids the triviality of
having factors equal to the identity.} In particular, given a
composition $\a$ and a sequence $\vec{i}=(i_2,i_3,\ldots)$ of
nonnegative integers (called the \bold{cycle index}), we wish to
determine the number $\Phur{\a}{i}$ of cycle factorizations of any
fixed permutation $\p \in \Class{\a}$ into exactly $i_2$ 2-cycles,
$i_3$ 3-cycles, \emph{etc}. Our results will be formulated in terms
of the generating series
\begin{equation}
\label{eq:hurwitzgs}
   \PshurGS{m}(\vec{x},\vec{q}, u) :=
    \sum_{n \geq 1}    \, \sum_{\vec{i} \,\geq\, \vec{0}}
    \sum_{\substack{\a \cmp n \\ \len{\a} =  m}}
       \Phur{\a}{\vec{i}} \, \vec{q}^{\vec{i}}\,
        \frac{x_1^{\a_1}}{\a_1} \cdots \frac{x_m^{\a_m}}{\a_m}
                                \frac{u^{r(\vec{i})}}{r(\vec{i})!}, \quad\qquad m \geq 1.
\end{equation}
Here, and throughout, $r(\vec{i}) := i_2+i_3 + \cdots$ denotes the
total number of factors  in any factorization counted by
$\Phur{\a}{\vec{i}}$, and $\vec{q} = (q_2, q_3, \ldots)$ is a vector
of indeterminates.

The structure of generic cycle factorizations is not well understood and, aside from explicit evaluations of
$\Phur{[n]}{\vec{i}}$ (see~\cite{gj-cactus, springer} and Theorem~\ref{thm:springer} of this paper), little
work has been done on their enumeration. However, a significant effort has been directed toward a natural
specialization of this problem, which is to count what we call \bold{$k$-cycle factorizations}. These are
cycle factorizations whose factors are all $k$-cycles for some fixed $k$.

The case $k=2$ (transposition factors) is particularly important geometrically. Counting \mbox{2-cycle}
factorizations of permutations is known as the \emph{Hurwitz problem}, and dates back to Hurwitz's original
investigations into the classification of almost simple ramified coverings of the sphere by the
sphere~\cite{hurwitz}.  The following formula, suggested but not completely proved by Hurwitz himself, gives
the number of \mbox{2-cycle} factorizations of any permutation of cycle type $(\a_1,\ldots,\a_m)$:
\begin{equation}
\label{eq:hurwitz}
    n^{m-3} (n+m-2)! \prod_{i=1}^m \frac{\a_i^{\a_i+1}}{\a_i!}.
\end{equation}
Although it has been extensively studied from various points of view, ranging from analytic to combinatorial
(see \cite{melou-schaeffer, goryunov-lando,gj-transitive}, for example), no purely bijective proof of this
striking enumerative formula is known.  Such a proof would be of tremendous interest as it could provide
further insight into the underlying geometry.  This is particularly true in light of a recent celebrated
result of Ekedahl, Lando, Shapiro, and Vainshtein~\cite{elsv} that identifies the enumeration of
\mbox{2-cycle} factorizations (\emph{i.e.} almost simple coverings) with the evaluation of certain Hodge
integrals, objects of great interest in the intersection theory of the moduli space of curves.
See~\cite{gj-gromov} for further details on these fascinating connections.  Indeed, the ultimate goal of our
approach to factorization problems is a full combinatorialization of Hodge integrals.

For arbitrary $k > 2$, counting \mbox{$k$-cycle} factorizations
 has not been as thoroughly examined and appears quite difficult.
Substantial progress was made in~\cite{gj-aspects}, where generating
series formulations for the number of such factorizations of
permutations with up to three cycles are given. In particular,
letting $\PkshurGS{m}{k}$ denote the series obtained by specializing
$\PshurGS{m}$ at $u=q_k=1$ and $q_i = 0$ for $i \neq k$, it is shown
there that $x\frac{d}{dx} \PkshurGS{1}{k}(x) = s(x)$ and
\begin{equation}
\label{eq:gjresult}
        \PkshurGS{2}{k}(x_1,x_2) =
        \log\pr{\frac{s(x_1)-s(x_2)}{x_1-x_2}}
            - \frac{s(x_1)^{k} - s(x_2)^k}{s(x_1)-s(x_2)},
\end{equation}
where $s \in \Q[[x]]$ (which depends on $k$) is the unique series
satisfying
\begin{equation}
\label{eq:seqn}
        s = x e^{s^{k-1}}.
\end{equation}
A more complicated expression for $\PkshurGS{m}{3}$ is also given in terms of $s$, but the calculations
necessary to evaluate $\PkshurGS{m}{k}$ for $k \geq 4$  become intractable.

\subsection{Equivalence up to Commutation of Disjoint Factors}
\label{ssec:equivalence}

There is a natural equivalence relation on cycle factorizations induced by permitting commutations of
disjoint factors. That is, we say two cycle factorizations are \bold{equivalent} if  one can be obtained from
the other by repeatedly exchanging adjacent factors that are disjoint in the sense that no symbol is moved by
both. For example, the following factorizations are equivalent:
\begin{align*}
    (3\,4\,7\,9) \cdot
    (6\,9\,7\,8) \cdot
    (1\,2) \cdot
    (2\,4\,5) \cdot (2\,6)
    \sim
    (1\,2) \cdot
    (3\,4\,7\,9) \cdot
    (2\,4\,5) \cdot
    (6\,9\,7\,8) \cdot
    (2\,6).
\end{align*}
Let $\IPhur{\a}{i}$ denote the number of inequivalent cycle
factorizations of $\p \in \Class{\a}$ with cycle index $\vec{i}$. We
shall study these numbers through the series
\begin{equation}
\label{eq:inequivgs}
   \IPshurGS{m}(\vec{x},\vec{q}, u) :=
    \sum_{n \geq 1}    \, \sum_{\vec{i} \,\geq\, \vec{0}}
    \sum_{\substack{\a \cmp n \\ \len{\a} =  m}}
       \IPhur{\a}{\vec{i}} \, \vec{q}^{\vec{i}}\,
        \frac{x_1^{\a_1}}{\a_1} \cdots \frac{x_m^{\a_m}}{\a_m}
                                u^{r(\vec{i})}, \quad\qquad m \geq
                                1.
\end{equation}
As before, let $\IPkshurGS{m}{k}(\vec{x})$ be the restricted series counting inequivalent \mbox{$k$-cycle}
factorizations obtained from $\IPshurGS{m}$ by setting $u=q_k=1$ and $q_i = 0$ for $i \neq k$ in
$\IPshurGS{m}$.

The problem of counting factorizations up to commutation can
apparently be traced back to Stanley, who originally posed it in the
context of \mbox{2-cycle} factorizations. The first result along
these lines came from Eidswick~\cite{eidswick} and
Longyear~\cite{longyear}, who proved (independently) that the number
of inequivalent \mbox{2-cycle} factorizations of the full cycle
$(1\,2\,\cdots\,n)$ is the generalized Catalan number
\begin{equation}
\label{eq:catalan}
    \frac{1}{2n-1}\binom{3n-3}{n-1}.
\end{equation}
Longyear's approach involved commutation of factorizations into a
canonical form.  This led to the following cubic functional equation
for the generating series $h(x) = \frac{d}{dx} \PkshurGS{1}{2}(x)$,
from which~\eqref{eq:catalan} is easily deduced:
\begin{equation}
\label{eq:hdefn}
    h(x) = 1 + xh(x)^3.
\end{equation}
Springer~\cite{springer} later generalized Longyear's argument to obtain an explicit formula for
$\IPhur{[n]}{i}$. His result is recovered here as Theorem~\ref{thm:ispringer}.  Also see~\cite{gj-macdonald}
for an alternative derivation of the number of inequivalent \mbox{$k$-cycle} factorizations of a full cycle.

More recently, Goulden, Jackson, and Latour~\cite{glj-inequivalent}
counted inequivalent \mbox{2-cycle} factorizations of any
permutation $\p \in \Class{[n,m]}$, proving that
\begin{equation}
\label{eq:igjresult}
    \PkshurGS{2}{2}(x_1,x_2)
    =
    \log\pr{1+x_1 x_2 h(x_1) h(x_2) \frac{h(x_1)-h(x_2)}{x_1-x_2}},
\end{equation}
where $h$ is defined by~\eqref{eq:hdefn}. Their method again relies on commutation to canonical form,  with
the additional aid of a clever combinatorial construction and an intricate inclusion-exclusion argument.

\subsection{Outline of the Paper and Statement of Results}
\label{ssec:outline}

The primary goal of this paper is to introduce a new technique in the enumeration of both cycle
factorizations and their equivalence classes under commutation.  We believe our approach to be of interest
for two principal reasons: First, it conveniently allows one to ignore the fine detail of factorizations
(\emph{i.e.} element-wise analysis) and focus on the grander structure, and second, it makes clearer the
structural parallels between the enumeration of factorizations and their equivalence classes. What follows is
a brief overview of the paper highlighting our main results.

In Section~\ref{sec:preliminaries} the reader is introduced to various constructs and conventions that are
used extensively throughout the paper. Our analysis of cycle factorizations then begins in
Section~\ref{sec:cyclefacts}, where we describe a graphical representation that allows~\eqref{eq:hurwitzgs}
to be viewed as a generating series for a special class of labelled planar maps. Cycle factorizations of full
cycles are then seen to correspond with particularly simple maps, namely \emph{cacti}, which are a natural
generalization of trees. This leads to the recovery of a known explicit formula for the number
$\Phur{[n]}{i}$ of cycle factorizations of a full cycle (see Theorem~\ref{thm:springer}).  It also marks our
first encounter with the series $\rt = \rt(x,\vec{q},u) \in \Q[\vec{q},u][[x]]$ defined as the unique
solution of the functional equation
\begin{equation}
\label{eq:rcdefn}
        \rt = x e^{u Q(\rt)},
\end{equation}
where $Q(z) \in \Q[\vec{q}][[z]]$ is given by
\begin{equation}
\label{eq:qdefn}
        Q(z) := \sum_{k \geq 2} q_k z^{k-1}.
\end{equation}
Clearly $\rt$ is a generalization of the series $s$ given by~\eqref{eq:seqn}, making it no surprise that it
plays a central role in our analysis. In particular, we shall present a graphical decomposition of maps
(called \emph{pruning}) that identifies an algebraic dependence of $\PshurGS{m}$ on $\rt$.  This is the
content of Theorem~\ref{thm:cactuspruning}, and is the centerpiece of our method. By exploiting the pruning
decomposition we deduce  following extension of~\eqref{eq:gjresult} to factorizations of arbitrary cycle
index.
\begin{thm}
\label{thm:mainthm1} Let $\rt$ and $Q$ be defined as above and, for
$i=1,2$, set $\rt_i = \rt(x_i, \vec{q},u)$. Then
\begin{align*}
    \PshurGS{2}(x_1,x_2,\vec{q},u) &= \log\pr{\frac{\rt_1-\rt_2}{x_1-x_2}}
                - u\frac{\rt_1 Q(\rt_1) - \rt_2 Q(\rt_2)}{\rt_1-\rt_2}.
\end{align*}
\qed
\end{thm}
In the end, our proof of Theorem~\ref{thm:mainthm1} still rests on an \emph{ad hoc} enumeration that we have not
been able to generalize.  Thus a formulation of $\PshurGS{m}$ for arbitrary $m$ remains out of reach.  We believe
that it will be more tedious than difficult to extend our methods to arrive at an expression for $\PshurGS{3}$,
but this has only yet been done in a special case.  We  comment further on these developments in
\S\ref{ssec:furtherresults}.

In Section~\ref{sec:inequiv} we turn to the enumeration of cycle factorizations up to the equivalence defined
in~\S\ref{ssec:equivalence}. By modifying our graphical representation of cycle factorizations to allow for
commutations of adjacent factors, we are led to Springer's formula~\cite{springer} for the number of inequivalent
cycle factorizations of a full cycle (Theorem~\ref{thm:ispringer}). We discover that the unique solution $\irt =
\irt(x,\vec{q},u) \in \Q[\vec{q},u][[x]]$ of the functional equation
\begin{equation}
\label{eq:ircdefn}
    \irt = 1 + u \irt Q(x\irt^2)
\end{equation}
plays a role directly analogous to that of $\rt$ in the enumeration of ordered cycle factorizations. Again we
develop a pruning decomposition (Theorem~\ref{thm:inequivpruning}) from which we deduce the following
generalization of~\eqref{eq:igjresult} to factorizations of arbitrary cycle index.

\begin{thm}
\label{thm:mainthm2} Let $\irt$ and $Q$ be defined as above
and, for $i=1,2$, set $\irt_i = \irt(x_i, \vec{q},u)$. Then
\begin{equation*}
    \IPshurGS{2}(x_1,x_2,\vec{q},u)
    =
    \log \pr{\frac{(x_1 \irt_1 - x_2 \irt_2)^2}{(x_1-x_2)(x_1\irt_1^2-x_2\irt_2^2)}}.
\end{equation*}
\qed
\end{thm}

Upon setting $u=q_2=1$ and $q_k = 0$ for $k \neq 2$, the defining equation~\eqref{eq:ircdefn} of $\irt$ transforms
to $\irt = 1 + x \irt^3$, thus identifying $\irt$ it with Longyear's series $h(x)$ (see~\eqref{eq:hdefn}). Under
these same specializations, it is easy to check that Theorem~\ref{thm:mainthm2} reduces to~\eqref{eq:igjresult}.

Again, we have been unable to extend Theorem~\ref{thm:mainthm2} to give a general expression for $\IPshurGS{m}$.
However, we have used our method to deduce a raw form of $\IPkshurGS{3}{2}$, thus extending the
Goulden-Jackson-Latour result in a different direction.  See \S\ref{ssec:ifurtherresults} for further comments on
this and related matters.

\section{Preliminaries}
\label{sec:preliminaries}

The definitions and notational conventions described in this section
are used extensively throughout the remainder.  We warn the reader
that some of these conventions are nonstandard.

\subsection{Cyclic Lists}
\label{ssec:cycliclists}

A \bold{cyclic list} is an equivalence class under the relation identifying finite sequences that are cyclic
shifts of one another. We write $\olist{a_1,\ldots,a_n}$ for the cyclic list with representative sequence
$(a_1,\ldots,a_n)$. Thus $\olist{a_1,a_2,a_3,a_4} = \olist{a_3,a_4,a_1,a_2}$.  Some liberties will be taken
with this notation; for instance, when considering the list $\olist{a_1,a_2,a_3,a_4}$, we adopt the
convention that $a_5$ is to be interpreted as $a_1$, while $a_0=a_4$, \emph{etc.}  A cyclic list
$L=\olist{a_1,\ldots,a_n}$ of real numbers is \bold{increasing} (respectively, \bold{nondecreasing}) if one
of its representative sequences is strictly increasing (nondecreasing). If $a_{i-1} \geq a_i$ then the pair
$(a_{i-1},a_i)$ is called a \bold{descent} of $L$.

\subsection{Maps and Polymaps}
\label{ssec:maps}

Recall that a \bold{planar map} (subsequently, a \bold{map}) is a 2-dimensional cellular complex whose polyhedron
is homeomorphic to the sphere. The 0-cells, 1-cells, and 2-cells of a map are its \bold{vertices}, \bold{edges},
and \bold{faces}, respectively. An \bold{isomorphism} of two maps is an orientation-preserving homeomorphism
between their polyhedra which sends $i$-cells to $i$-cells and preserves incidence. We always consider isomorphic
maps to be indistinguishable.  If a map is labelled, then this definition of isomorphism is amplified to preserve
all labels.

The \bold{boundary walk} of a face $f$ of a map is the cyclic sequence $\olist{(v_0,e_0),\ldots,(v_{k},e_{k})}$ of
alternating vertices and edges listed in order as they are encountered along a counterclockwise traversal of the
boundary of $f$. (Counterclockwise here means that $f$ is always kept to the left of the line of traversal.) A
subsequence $(e_{i-1},v_i,e_i)$ consisting of two consecutive edges and their common incident vertex is called a
\bold{corner} of $f$.

\begin{figure}[t]
    \begin{center}
    \includegraphics[width=0.65\textwidth]{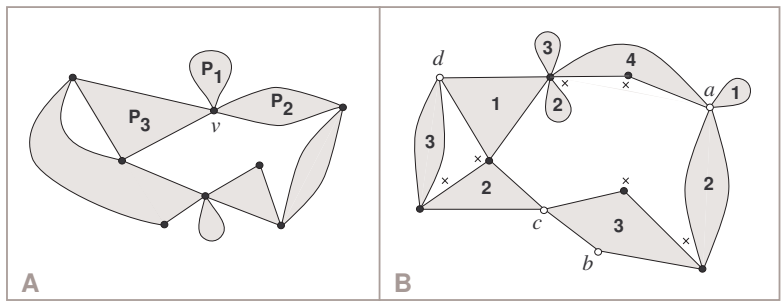}
    \caption{(A) A polymap with 8 polygons and 3 faces. (B) The descent structure of a polymap}
    \label{fig:polymap}
    \end{center}
\end{figure}

A map is \bold{2-coloured} if its faces have been painted black and white so every edge is incident with both
a black face and a white face (so no two similarly coloured faces are adjacent). We are interested in the
special class 2-coloured maps for which the boundary walk of every black face is a cycle (\emph{i.e.}
contains no repeated vertices or edges). We call these \bold{polymaps}, and make the following supporting
definitions:

\begin{itemize}
\setlength{\itemsep}{0pt}
\item   The black faces of a polymap are called \bold{polygons}. An \bold{$m$-gon} is a black face of
        degree $m$.
\item   The white faces of a polymap are referred to simply as \bold{faces}.

\item   A \bold{corner} of a polymap always refers to a corner of a (white) face.

\item   The \bold{rotator} of a vertex $v$ in a polymap is the unique cyclic list of polygons encountered along a
        clockwise tour of small radius about $v$.
\end{itemize}
For example, a polymap with 9 polygons and 3 faces is illustrated in
Figure~\ref{fig:polymap}A. The rotator of vertex $v$ is
$\olist{p_1,p_2,p_3}$.

\subsection{Descent Structure}
\label{ssec:descents}

Let $\amap$ be a polymap whose polygons are labelled with integers. We shall  always regard the edges of
$\amap$ as being labelled, with each edge inheriting its label from the unique polygon that it borders. Let
$f$ be a face of the $\amap$ with boundary walk $\olist{(v_0,e_0),\ldots,(v_k,e_k)}$. If $e_{i-1} \geq e_i$,
then the pair $(e_{i-1},e_i)$ a \bold{descent} of $f$.  We say $v_i$ is \bold{at} this descent, and the
corner $(e_{i-1},v_i,e_i)$  is called a \bold{descent corner}.  The \bold{descent set} of $f$ is the set
$\{v_i \,:\, e_{i-1} \geq e_i\}$ of all vertices at descents of $f$. The cyclic list obtained by listing the
vertices of this set in the order in which they appear along the boundary walk of $f$ is called the
\bold{descent cycle} of $f$.

\begin{exmp} Consider the outer face $f$ of the polymap in Figure~\ref{fig:polymap}B.  The cyclic list of edge labels encountered along
its boundary walk is $\olist{1,2,3,3,2,3,1,3,4}$, so $f$ has 4
descents, namely $(3,3), (3,2), (3,1)$, $(4,1)$. The hollow vertices
$a,b,c$, and $d$ are at these descents, so the descent set of $f$ is
$\{a,b,c,d\}$ and its descent cycle is $\olist{a,b,c,d}$. The
descent corners of the other faces are marked with crosses. \qed
\end{exmp}

\subsection{Constellations}
\label{ssec:constellations}

\begin{figure}[t]
    \begin{center}
    \includegraphics[width=.35\textwidth]{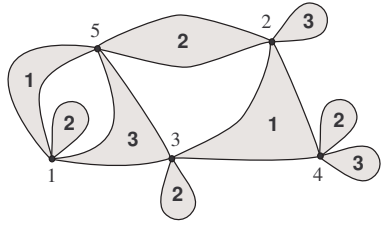}
    \caption{A 3-constellation on 5 vertices.} \label{fig:constellation}
    \end{center}
\end{figure}

An \bold{$r$-constellation} on $n$ vertices is a polymap whose vertices are distinctly labelled
$1,2,\ldots,n$, and whose polygons are labelled $1,\ldots,r$ such that the rotator of every vertex is
$\olist{1,2,\ldots,r}$. Figure~\ref{fig:constellation} illustrates a $3$-constellation on $5$ vertices.

Our interest in this special class of polymaps stems from a well
known connection between maps and factorizations. A complete
description of this correspondence can be found
in~\cite{melou-schaeffer} (albeit in dual form\footnote{We have
borrowed the term ``constellation'' from~\cite{melou-schaeffer},
though we caution that the term is used there for an object that is
dual to those considered here.}), so here it will suffice to outline
a version particularly suited to our needs.

\begin{prop}
\label{prop:mainbijection} Minimal transitive factorizations in
$\Sym{n}$ with $r$ factors are in bijection with
\mbox{$r$-constellations} on $n$ vertices. In particular,
factorizations of $\p$ correspond with constellations whose decent
cycles are exactly the cycles of $\p$.
\end{prop}

\begin{proof}[Sketch proof:]
Let $\amap[C]$ be an $r$-constellation on $n$ vertices. From $\amap[C]$ define permutations
$\s_1,\ldots,\s_r$ as follows:  To form $\s_i$, first obtain a collection of disjoint cycles by listing the
vertices of each polygon labelled $i$ as they appear in clockwise order around its perimeter, and then let
$\s_i$ be the product of these cycles. (See Example~\ref{exmp:constellation}, below.)

Now set $\p=\s_r \cdots \s_1$. Regarding $\p$ as the sequential product of $\s_1,\ldots,\s_r$, observe that
it acts on a vertex $v$ to move it clockwise around the polygons of $\amap[C]$ following edges labelled
$1,2,\ldots,r,$ in turn. That is, $v$ starts at the descent corner $(r,v,1)$ of some face and is moved by
$\pi$ along the boundary walk until it lands at the next descent corner $(r,v',1)$ of that face. Since
$\p(v)=v'$, the cycles of $\p$ are the descent cycles of $\amap[C]$.

Set $F=(\s_r,\ldots,\s_1)$. We claim $\amap[C] \leftrightarrow F$ is the desired correspondence.  Clearly $F$
is a factorization of $\p$, and it is transitive because $\amap[C]$ is connected. Note that $\amap[C]$ has
$n$ vertices and  $nr$ edges. If $\p \in \Class{\a}$ and $\s_i \in \Class{\b_i}$ for all $i$ then $\amap[C]$
has $\len{\a} + \sum_i \len{\b_i}$ total faces (white faces $+$ polygons). Since $\amap[C]$ is planar,
Euler's polyhedral formula ($\#\text{vertices}-\#\text{edges}+\#\text{faces}=2$) gives equality
in~\eqref{eq:genusdefn}, so $F$ is minimal transitive.
\end{proof}

\begin{exmp} \label{exmp:constellation} The $3$-constellation of Figure~\ref{fig:constellation} corresponds to the factorization
$F = (\s_3, \s_2, \s_1)$ of $\p=(1\,2\,4)(3)(5)$, where $\s_1 = (1\,5)(2\,4\,3)$, $\s_2=(1)(2\,5)(3)(4)$, and
$\s_3=(1\,5\,3)(2)(4)$. Note that vertices $1,2$ and $4$ are at descents of the outer face, and the cycle
$(1\,2\,4)$ of $\p$ coincides with the descent cycle of that face.\qed
\end{exmp}

\section{Cycle Factorizations}
\label{sec:cyclefacts}

In this section we shall introduce a convenient graphical representation of cycle factorizations as decorated
polymaps, and proceed to give a decomposition of these polymaps that can simplify their enumeration.  The
section culminates in a proof of Theorem~\ref{thm:mainthm1}.

\subsection{Proper and $\a$-Proper Polymaps}
\label{ssec:properpolymaps}

\begin{figure}[t]
    \begin{center}
    \includegraphics[width=.85\textwidth]{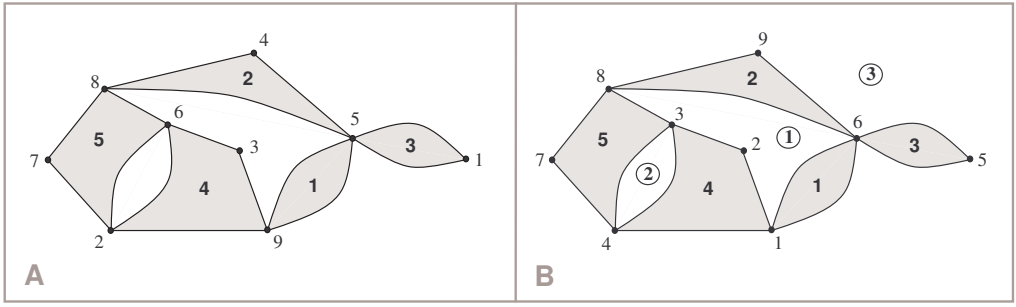}
    \caption{(A) The polymap of a cycle factorization.  (B) A $(3,1,5)$-proper polymap. } \label{fig:properpolymap}
    \end{center}
\end{figure}

Let $F=(\s_r,\ldots,\s_1)$ be a cycle factorization and let $\amap[C]$ be the constellation corresponding to $F$
through the bijection of Proposition~\ref{prop:mainbijection}.  Then exactly $r$ of the polygons of $\amap[C]$ are
of degree two or more, while the remainder are 1-gons bounded by loops. Let $\amap_F$ be the polymap obtained by
removing all 1-gons.  For example, Figure~\ref{fig:properpolymap}A
shows $\amap_F$ when $F$ is the factorization
\begin{equation}
\label{eq:cycexmp2}
    (1\,5\,7\,8\,4)(2)(3\,9\,6)
    =
    (2\,7\,8\,6) \cdot (2\,6\,3\,9)\cdot(1\,5)\cdot(4\,5\,8)\cdot(5\,9).
\end{equation}
Notice that all rotators of $\amap_F$ are increasing, since every rotator of $\amap[C]$ is
$\olist{1,\ldots,r}$. In fact, the removal of loops does not affect descent structure, in the sense that
$\amap_F$ has precisely the same descent cycles as $\amap[C]$. Moreover, $\amap[C]$ is easily recovered from
$\amap_F$ by attaching loops in the unique manner which makes each rotator $\olist{1,\ldots,r}$.

We say  a loopless polymap is \bold{proper} if its polygons are labelled with distinct integers so that all
rotators are increasing. (The polygon labels are assumed to be $1,\ldots,r$, for some $r$, unless otherwise
specified.) Let the \bold{poly-index} of a loopless polymap be the vector $\vec{i}=(i_2,i_3,\ldots)$, where
$i_k$ is the number of $k$-gons it contains. With these definitions, we have:

\begin{prop}
\label{prop:mapping} The mapping $F \mapsto \amap_F$ is a bijection between cycle factorizations in $\Sym{n}$ and
proper polymaps with vertices labelled $1,\ldots,n$. In particular, if $F$ is a factorization of $\p$ of cycle
index $\vec{i}$, then $\amap_F$ is a polymap of poly-index $\vec{i}$ whose descent cycles are precisely the cycles
of $\p$. \qed
\end{prop}

\newcommand{\Pai}{\mathcal{P}_{\!\a}(\vec{i})}

Let $\a=(\a_1,\ldots,\a_m)$ be a composition, and let $\Pai$ be the set of
vertex-labelled proper polymaps of poly-index $\vec{i}$ whose
descent sets are $\Dset{\a}{1},\ldots,\Dset{\a}{m}$, where
\begin{align}
\label{eq:canonicalsets}
    \Dset{\a}{j} := \{k \in \N \,:\, \a_1+\cdots+\a_{j-1} < k \leq \a_1+\cdots+\a_j\}.
\end{align}
Under the mapping $F \mapsto \amap_F$, members of $\Pai$ correspond with factorizations of permutations whose
orbits are $\Dset{\a}{1},\ldots,\Dset{\a}{m}$.  Since there are $\prod_j (\a_j-1)!$ such permutations, each
admitting the same number of factorizations of given cycle index, we conclude that $|\Pai| = \Phur{\a}{i}
\cdot \prod_j (\a_j-1)!$.

Let us say a proper polymap is \bold{$\a$-proper} if its faces are labelled $1,\ldots,m$ so that face $j$ has
$\a_j$ descents, for $1 \leq j \leq m$. (See Figure~\ref{fig:properpolymap}B, for example.)  Notice that any
$\a$-proper polymap of poly-index $\vec{i}$ can be transformed into a member of $\Pai$ by first labelling its
vertices in any of the $\smash{\prod_j \a_j!}$ ways that make $\Dset{\a}{j}$ the descent set of face $j$ for
$j=1,\ldots,m$, and then stripping face labels.  In fact, no members of $\Pai$ are duplicated in this process
when $\len{\a} \geq 2$, since the face and polygon labels of $\a$-proper polymaps preclude nontrivial
automorphisms.\footnote{This is not so when $\len{\a}=1$, since a polymap composed of a single polygon has
rotational symmetry.} So for $\len{\a} \geq 2$ we have $|\Pai| = \Pshur{\a}{i} \cdot \prod_j \a_j!$, where
$\Pshur{\a}{i}$ is the number of $\a$-proper polymaps of poly-index $\vec{i}$.

Comparing the expressions above for $|\Pai|$ gives $\Phur{\a}{i} =
\a_1\cdots\a_m \Pshur{\a}{i}$ when $\len{\a} \geq 2$.
Thus~\eqref{eq:hurwitzgs} becomes
\begin{equation}
\label{eq:propergs}
    \PshurGS{m}(\vec{x},\vec{q}, u) =
    \sum_{n \geq 1}    \, \sum_{\vec{i} \,\geq\, \vec{0}}
    \sum_{\substack{\a \cmp n \\ \len{\a} =  m}}
       \Pshur{\a}{i} \, \vec{q}^{\vec{i}}\,
        \vec{x}^{\pmb{\a}}
        \frac{u^{r(\vec{i})}}{r(\vec{i})!},\qquad\text{for $m \geq
        2$}.
\end{equation}
That is, $\PshurGS{m}(\vec{x},\vec{q},u)$ is the generating series
for $\a$-proper polymaps, where $u$ is an exponential marker for labelled polygons, $x_j$
marks descents of face $j$ (for $1 \leq j \leq m$), and $q_k$
records the number of $k$-gons (for $k \geq 2$).

\subsection{Proper Cacti and  Factorizations of Full Cycles}
\label{ssec:fullcycles}

\newcommand{\hrt}{\bar{\rt}}

A \bold{cactus} is a polymap with only one face. Hence Proposition~\ref{prop:mapping} implies $\Phur{[n]}{i}$
is the number of proper cacti of poly-index $\vec{i}$ whose sole faces have descent cycle
$(1\,2\,\cdots\,n)$. However, observe that the fixed descent cycle of such a cactus forces all vertex labels
once the location of one is known. Thus $\Phur{[n]}{i}$ is the number of \emph{vertex-rooted} proper cacti of
poly-index $\vec{i}$. (That is, the root marks the location of a canonical label, say 1.)

Let $\rt=\rt(x,\vec{q},u)$ be the generating series for
vertex-rooted proper cacti with respect to total vertices (marked by
$x$), labelled polygons (marked by $u$), and poly-index (marked by
$\vec{q}$). That is,
\begin{equation} \label{eq:rcexpansion}
   \rt
        :=\sum_{n \geq 1} \sum_{\vec{i} \geq \vec{0}} \Phur{[n]}{\vec{i}} \,\vec{q}^{\vec{i}}
                            {x^n} \frac{u^{r(\vec{i})}}{r(\vec{i})!}
        = x\frac{d}{dx} \PshurGS{1}(x,\vec{q},u).
\end{equation}
Though this may appear to be in conflict with our earlier definition~\eqref{eq:rcdefn} of $\rt$, we now give
a combinatorial decomposition of cacti to show that the two definitions coincide.

\begin{figure}[t]
    \begin{center}
    \includegraphics[width=\textwidth]{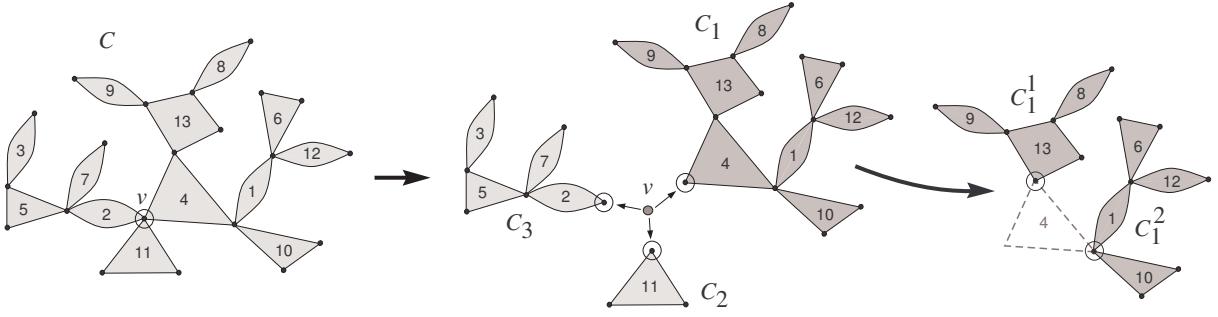}
    \caption{Decomposition of a rooted cactus.} \label{fig:polycactusdecomp1}
    \end{center}
\end{figure}

Suppose that the root $v$ of a rooted proper cactus $C$ is incident
with $m$ polygons. Detaching these polygons from $v$ results in a
collection $\{C_1,\ldots,C_m\}$ of rooted proper cacti, where the
root of each $C_i$ is incident with only one polygon. (See
Figure~\ref{fig:polycactusdecomp1}.)
The ordering of $C_1,\ldots,C_m$ around $v$ need not be recorded, as it can be deduced from the increasing rotator
condition. Thus
\begin{equation}
\label{eq:hrt}
    \rt = x \sum_{m \geq 0} \frac{\hrt^m}{m!} = x e^{\hrt},
\end{equation}
where the series $\hrt=\hrt(x,\vec{q},u)$ counts rooted proper cacti
such as $C_i$ in which $x$ marks only \emph{nonroot} vertices. But
if the root of $C_i$ is incident with a $k$-gon, then removal of
this polygon leaves a $(k-1)$-tuple $C^{1}_{i},\ldots,C^{k-1}_i$ of
rooted proper cacti. (See Figure~\ref{fig:polycactusdecomp1},
right-hand side.) This accounts for a contribution $u q_k \rt^{k-1}$
to $\hrt$, and summing over $k$ yields $\hrt=uQ(\rt)$, where $Q$ is
defined as in~\eqref{eq:qdefn}. Thus~\eqref{eq:hrt} gives $\rt =
xe^{uQ(\rt)}$, in agreement with~\eqref{eq:rcdefn}.

The following result was originally proved bijectively by Springer~\cite{springer}.  In fact, the
correspondence between factorizations and cacti employed here was also developed in~\cite{springer}, but from
a less general point of view.

\begin{thm}[Springer~\cite{springer}]
\label{thm:springer} Let $\vec{i}=(i_2,i_3,\ldots)$ be a sequence of nonnegative integers and set
$r=r(\vec{i})=i_1+i_2+\cdots$. Then the number of cycle factorizations of $(1\,2\,\cdots\,n)$ with cycle index
$\vec{i}$ is
$$
\Phur{[n]}{\vec{i}} = \frac{n^{r-1}\,r!}{\prod_{k \geq 2} i_k!}
$$
in the case that $n + r - 1 = \sum_{k \geq 2} ki_k$, and zero  otherwise. \qed
\end{thm}

\begin{proof}
From~\eqref{eq:rcexpansion} we have $\Phur{[n]}{\vec{i}} = r!\,[x^n
\vec{q}^{\vec{i}} u^r]\,\rt(x,\vec{q},u)$. This coefficient can be
obtained from $~\eqref{eq:rcdefn}$ through a straightforward
application of Lagrange inversion~\cite{goulden-jackson}. An
alternative proof that relies on a Pr\"ufer-like encoding of cacti
can be found in~\cite{springer}.
\end{proof}


\subsection{Smooth Polymaps, Cores, and Branches}
\label{ssec:smoothpolymaps}

Given~\eqref{eq:propergs}, we now wish to count $\a$-proper polymaps for general $\a$ with $\len{\a} \geq 2$.
To do so we invoke a technique we call \emph{pruning}, whereby polymaps are simplified through the removal of
cacti. In this section we lay the groundwork for this approach, which is then executed in the next.

\begin{figure}[t]
\begin{center}
\psfrag{frag:m}{$\amap$} \psfrag{frag:b}{$B$} \psfrag{frag:mc}{$\acore$}
\includegraphics[width=.85\textwidth]{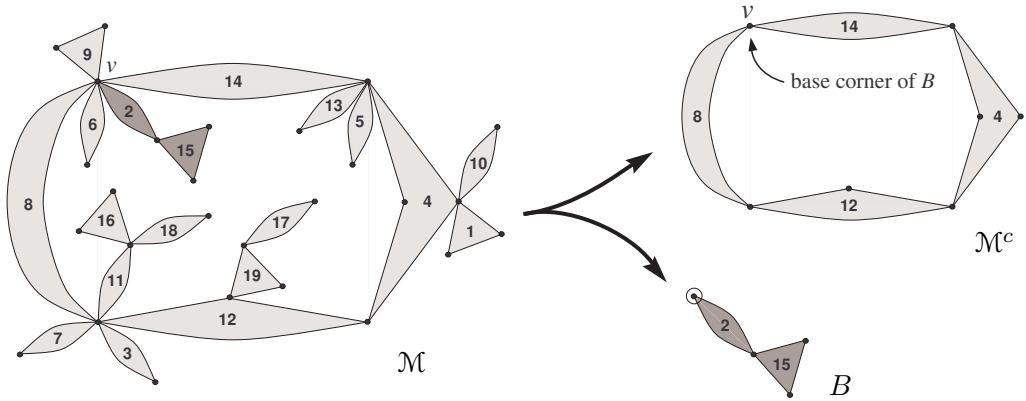}
\caption{A proper polymap, its core, and one of its branches.} \label{fig:polymapcore}
\end{center}
\end{figure}

A \bold{leaf} of a polymap is a polygon that shares exactly one vertex with another polygon,  and a polymap
is \bold{smooth} if it does not have any leaves. If $\amap$ is a polymap with at least two faces, then
iteratively removing its leaves results in a unique smooth polymap called the \bold{core} of $\amap$ and
denoted by $\acore$. Labels of $\acore$ are inherited from $\amap$ in the obvious way. See
Figure~\ref{fig:polymapcore}, for example. If $p$ is a polygon of $\amap$ sharing only one vertex $v$ with
$\acore$, then separating $p$ from $v$ results in two components, one of which is a rooted cactus $B$ whose
root vertex is incident only with $p$. We call $B$ a \bold{branch} of $\amap$, and refer to $p$ as its
\bold{root polygon}.  The \bold{base corner} of $B$ is the corner of $\acore$ in which $p$ was attached. (See
Figure~\ref{fig:polymapcore}.)

Consider now the case when $\amap$ is proper, with at least two
faces. Let $B$ and $p$ be as described above, and let $f$ be the
face of $\acore$ in which $B$ was attached.  Let
$\olist{(v_0,e_0),\ldots,(v_k,e_k)}$ be the boundary walk of $f$,
where the indexing has been chosen so that the sequence
$(e_0,\ldots,e_k)$ is as small as possible in lexicographic order.
It is not difficult to see that this condition specifies
$e_0,\ldots,e_k$ uniquely, and hence also determines a unique value
of $b$, with $0 \leq b \leq k$, such that $(e_{b-1},v_b,e_b)$ is the
base corner of $B$.

The key observation here is that the increasing rotator condition on $\amap$ implies $\olist{e_{b-1},p,e_b}$
is nondecreasing, and that this puts strong restrictions on the possible values of $b$, as is demonstrated by
the following lemma.

\begin{lem}
\label{lem:cyclic} Let $L = \olist{a_0,\ldots,a_k}$ be a cyclic list of real numbers with $d$ descents. If $a \in
\R$ is not in $L$, then there are exactly $d$ values of $i$ with $0 \leq i \leq k$ such that
$\olist{a_{i-1},a,a_{i}}$ is nondecreasing.
\end{lem}

\begin{figure}[t]
    \begin{center}
    \includegraphics[width=.35\textwidth]{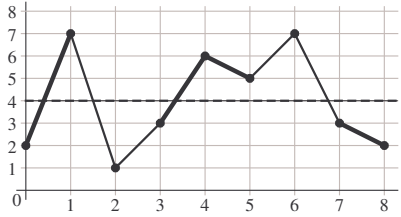}
    \caption{The proof of Lemma~\ref{lem:cyclic}, with $L=\olist{2,7,1,3,6,5,7,3}$ and $a=4$.} \label{fig:cyclelemma}
    \end{center}
\end{figure}

\begin{proof}
Let $\mathcal{P}$ be the polygonal path in the plane connecting the
points $(0,a_0),\ldots,(k,a_k),(k+1,a_0)$, in that order. Let $s_i$
be the $i$-th step of $\mathcal{P}$.  Call $s_i$ an \emph{up step}
if $a_i > a_{i-1}$ and a \emph{down step} otherwise.  Note that
$\olist{a_{i-1},a,a_i}$ is nondecreasing if and only if either
$a_{i-1} < a< a_i$, or $a > a_{i-1} \geq a_i$, or $a_{i-1} \geq a_i
> a$. Plainly, one of these conditions holds if and only if either
(1) $s_i$ is an up step which the line $y = a$ crosses, or (2) $s_i$
is a down step which this line misses. Since $\mathcal{P}$ begins
and ends at the same $y$-coordinate, the numbers of up steps and
down steps crossed by $y=a$ must be equal. Thus the number of
indices $i$ for which (1) or (2) is satisfied is equal to number of
down steps of $P$. Since down steps reflect descents of $L$, this
completes the proof. See Figure~\ref{fig:cyclelemma} for an
illustration; the dashed line is $a=4$, and steps for which
$\olist{a_{i-1},a,a_i}$ is nondecreasing have been thickened.
\end{proof}

In particular, the lemma implies $\olist{e_{s-1},p,e_s}$ is
increasing for exactly $d$ values of $s$ in the range $0 \leq s \leq
k$, where $d$ is the number of descents of face $f$.  Let these
possible $s$-values be $s_1 < \cdots < s_d$. Then since
$\olist{e_{b-1},p,e_b}$ is nondecreasing, we have $b=s_i$ for a
unique  $i \in \{1,\ldots,d\}$. We call $i$ the \bold{index} of the
branch $B$.

\begin{exmp}  Consider branch $B$ of $\amap$ shown in
Figure~\ref{fig:polymapcore}, and let $f$ be the face of $\acore$ in
which $B$ was attached. The boundary walk of $f$ is
$\olist{(v_0,e_0),\ldots,(v_5,e_5)}$, where we choose
$(e_0,\ldots,e_5)=(4,4,14,8,12,12)$ to be lexicographically minimal.

The base corner of $B$ is $(e_2,v_3,e_3)=(14,v_3,8)$ and the polygon
attached to its root has label 2, so in the notation used above we
have $b=3$ and $p=2$. Note that $f$ has $d=4$ descents, namely
$(14,8)$, $(12,12)$, $(12,4)$, and $(4,4)$.  In accordance with the
lemma, $\olist{e_{s-1}, p, e_s}$ is nondecreasing for exactly $d=4$
values of $s \in \{0,\ldots,5\}$, namely $s_1=1$, $s_2=3$, $s_3=4$,
and $s_4=5$. Since $b=s_2$, the index of $B$ is $i=2$. \qed
\end{exmp}

To reiterate, Lemma~\ref{lem:cyclic} guarantees that a given branch can be attached to a face with $d$
descents in \emph{exactly} $d$ positions so as to maintain increasing rotators, and the index of a branch
records which of these positions it occupies.

\subsection{Pruning Cacti}
\label{ssec:pruningcacti}

Let $\Psmhur{\th}{i}$ denote the number of smooth $\th$-proper
polymaps with poly-index $\vec{i}$. For $m \geq 2$ we define the
following smooth polymap analogue of the series $\PshurGS{m}$:
\begin{equation}
\label{eq:smoothgs}
    \PsmhurGS{m}(\vec{x},\vec{q},u) :=
    \sum_{n \geq 1}    \, \sum_{\vec{i} \,\geq\, \vec{0}}
    \sum_{\substack{\th \cmp n \\ \len{\th} =  m}}
            \Psmhur{\th}{i} \, \vec{q}^{\vec{i}} \,
                                \vec{x}^{\pmb{\th}}
                                \frac{u^{r(\vec{i})}}{r(\vec{i})!}.
\end{equation}

The following theorem is the centrepiece of our approach to counting cycle factorizations.  At first glance
it appears to be the rather transparent statement that arbitrary polymaps can be viewed as the composition of
smooth polymaps with cacti. However, it is worth emphasizing that we are concerned with the decomposition of
\emph{proper} polymaps, and therefore  must maintain control of descent structure through the pruning
process.  It is Lemma~\ref{lem:cyclic}, and the resulting ``index of a branch'', that enables us to deal with
this subtlety.

\begin{thm}
\label{thm:cactuspruning}  Fix $m \geq 2$ and set $\rt_i =
\rt(x_i,\vec{q},u)$ for $i=1,\ldots,m$, where $\rt$ is given
by~\eqref{eq:rcdefn}. Then
\begin{equation*}
\label{eq:cactipruning}
    \PshurGS{m}(\vec{x}, \vec{q},  u) =
    \PsmhurGS{m}(\vec{\rt}, \vec{q}, u),
\end{equation*}
where $\vec{x} = (x_1,\ldots,x_m)$ and $\vec{\rt} = (\rt_1,\ldots,\rt_m)$.
\end{thm}

\begin{proof}
Let $\a$ be an $m$-part composition, let $\amap$ be an $\a$-proper
polymap. Suppose $\acore$ has $\th_j$ descents in face $j$, for $1
\leq j \leq m$. Then for all $1 \leq j \leq m$ and $1 \leq i \leq
\th_j$, let $\aset{B}^j_i$ be the set of all branches of index $i$
in face $f$ of $\amap$. Assemble all branches of $\aset{B}^j_i$ into
a single rooted proper cactus $C^j_i$ by identifying their root
vertices, and  let $\aset{O}_j$ be the ordered forest
$(C^j_1,\ldots,C^j_{\th_j})$. See Figure~\ref{fig:pruningcacti} for
an example of these constructions.
\begin{figure}[t]
    \begin{center}
    \psfrag{frag:m}{$\amap$} \psfrag{frag:mc}{$\acore$} \psfrag{frag:f1}{$\aset{O}_1$} \psfrag{frag:f2}{$\aset{O}_2$}
    \psfrag{frag:b11}{$\aset{B}^1_1$} \psfrag{frag:b12}{$\aset{B}^1_2$} \psfrag{frag:b13}{$\aset{B}^1_3$}
    \psfrag{frag:b14}{$\aset{B}^1_4$} \psfrag{frag:b21}{$\aset{B}^2_1$} \psfrag{frag:b22}{$\aset{B}^2_3$}
    \psfrag{frag:b23}{$\aset{B}^2_3$}
    \includegraphics[width=\textwidth]{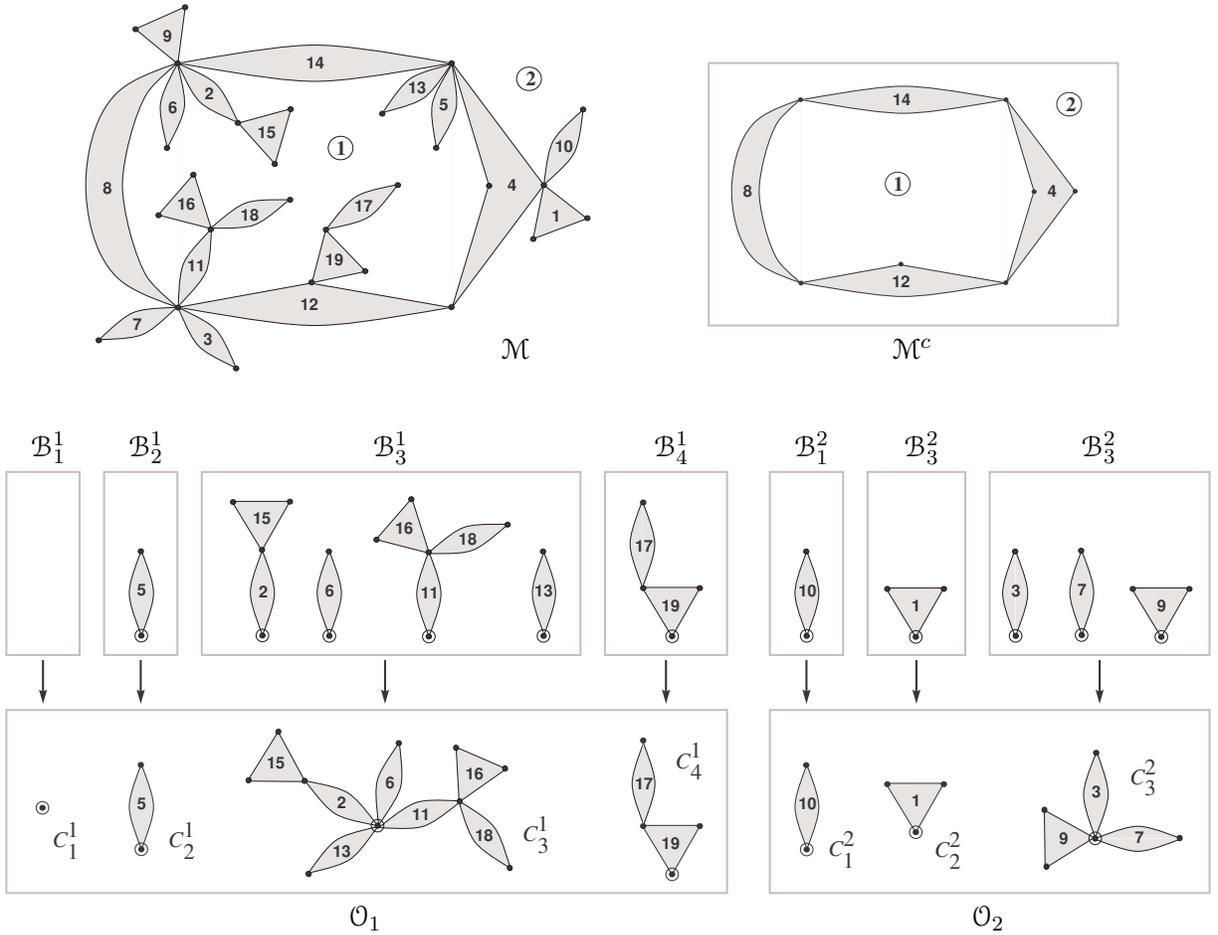}
    \caption{Pruning cacti from an $\a$-proper polymap.} \label{fig:pruningcacti}
    \end{center}
\end{figure}

We claim that the mapping $\Phi \,:\, \amap
\mapsto (\aset{O}_1,\ldots,\aset{O}_m, \acore)$ is a polygon-preserving
bijection between $\a$-proper polymaps and tuples
$(\aset{F}_1,\ldots,\aset{F}_m,
\amap[S] )$ satisfying the following properties:
\begin{itemize}
\setlength{\itemsep}{0pt}
\item[(a)]  $\aset{F}_j$ is an ordered forest containing $\th_j$ rooted proper cacti with a total of $\a_j$  vertices. 

\item[(b)]  $\amap[S]$ is a smooth $\th$-proper polymap, where $\th=(\th_1,\ldots,\th_m)$.

\item[(c)]  The polygon labels of $\amap[S]$ and $\aset{F}_1,\ldots,\aset{F}_m$ together
            partition $\{1,\ldots,r\}$.
\end{itemize}
The fact that $\Phi(\amap)$ is indeed a tuple of this type follows from two simple observations: (1) a vertex is
at a descent of face $f$ of $\amap$ if and only if it is at a descent of face $f$ of $\acore$, and (2) every
non-root vertex of a branch of $\amap$ is at a descent of the face in which it lies.

The injectivity of $\Phi$ is immediate because two $\a$-proper polymaps are isomorphic if and only if their
cores are isomorphic and all branches in corresponding faces agree.  To prove $\Phi$ is also surjective, let
$(\aset{F}_1,\ldots,\aset{F}_m,\amap[S])$ satisfy (a) through (c), where $\smash{\aset{F}_j =
(C^j_1,\ldots,C^j_{\th_j})}$ for $1 \leq j \leq m$. For fixed $j$ and $i$ with $1 \leq i \leq \th_j$, the
rooted cactus $C^j_i$ can unambiguously be viewed as a collection of branches whose roots have been
identified. Take any such branch, $B$, and let $p$ be the label of its root polygon. Now let
$\olist{(v_0,e_0),\ldots,(v_k,e_k)}$ be the boundary walk of face $j$ of $\amap[S]$, with $(e_0,\ldots,e_k)$
lexicographically minimal.  Then $\olist{e_0,\ldots,e_k}$ has $\th_j$ descents according to (b). Since (c)
ensures that $p$ is distinct from $e_0,\ldots,e_k$, Lemma~\ref{lem:cyclic} implies $\olist{e_{s-1},p,e_s}$ is
nondecreasing for exactly $\th_j$ values of $s$, say $0 \leq s_1 < \cdots < s_{\th_j} \leq k$.  Attach $B$ to
face $j$ of $\amap[S]$ at corner $(e_{s_i-1},v_{s_i},e_{s_i})$, doing so in the unique manner that leaves the
rotator of $v_{s_i}$ increasing. Repeat this process for all $i$, $j$ and $B$ to iteratively build a polymap
$\amap$. Note that the order in which branches are attached is immaterial, so $\amap$ is well defined. The
fact that $\amap$ is $\a$-proper follows from conditions (a)--(c) and observations (1) and (2) made above.
Since $\Phi(\amap)=(\aset{F}_1,\ldots,\aset{F}_m, \amap[S])$, by construction, we conclude that $\Phi$ is
surjective.

This bijection shows $\Pshur{\a}{i}$ to be the number of tuples $(\aset{F}_1,\ldots,\aset{F}_m,\amap[S])$
satisfying (a)--(c) and with poly-index $\vec{i}$.  The result now follows by comparing~\eqref{eq:propergs}
and~\eqref{eq:smoothgs}, and recalling from~\S\ref{ssec:fullcycles} that $\rt$ is the generating series for
rooted proper cacti.
\end{proof}

\subsection{Factorizations of Permutations with Two Cycles}
\label{ssec:twofacepolymaps}

We now conclude this section with a proof of Theorem~\ref{thm:mainthm1}. The proof relies on
Theorem~\ref{thm:cactuspruning} and the following well known result, which can be found
in~\cite{goulden-jackson}. Recall that a \emph{circular permutation} of $\{1,\ldots,n\}$ is a cyclic list
$\olist{a_1,\ldots,a_{n}}$ such that $\{a_1,\ldots,a_{n}\} = \{1,\ldots,n\}$.

\begin{lem}
\label{lem:circperm} There are
$$
        n!\,[x^{n-d} y^d]\,\log\pr{\frac{x-y}{xe^y-ye^x}}
$$
circular permutations of $\{1,\ldots,n\}$ having exactly $d$
descents. \qed
\end{lem}

\vspace*{1.5mm}

\begin{figure}[t]
    \begin{center}
    \includegraphics[width=.25\textwidth]{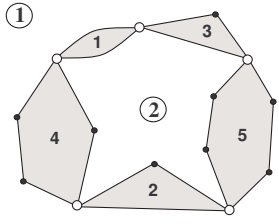}
    \caption{A smooth proper polymap.} \label{fig:2facepolymap}
    \end{center}
\end{figure}

\noindent\textbf{Proof of Theorem~\ref{thm:mainthm1}:} We first
determine the series $\PsmhurGS{2}(z_1,z_2,\vec{q},u)$ counting
smooth $\a$-proper polymaps with $\len{\a}=2$. Let $\amap[S]$ be
such a polymap, noting that this means $\amap[S]$ is simply a closed
chain of polygons, each incident with exactly two others. (See
Figure~\ref{fig:2facepolymap}.  Vertex labels have been suppressed
for clarity.)
Let $\olist{l_1,\ldots,l_r}$ be the cyclic list of distinct polygon labels encountered along the boundary walk of
face 1. Set $\g_i=(j_1,j_2)$ if the polygon with label $l_i$ is a $(j_1+j_2)$-gon that has $j_1-1$ vertices
incident only with face $1$ and $j_2-1$ incident only with face $2$. Then $\amap[S]$ is fully specified by
$\olist{l_1,\g_1,\ldots,l_r,\g_r}$. For example, the polymap of Figure~\ref{fig:2facepolymap} corresponds with the
cyclic list
$$
\olist{1,(1,1),3,(2,1),5,(3,3),2,(1,2),4,(3,2)}.
$$
Let us say a vertex incident with only one face is \emph{internal} to that face; the remaining vertices are
\emph{extremal}. Extremal vertices are hollow in Figure~\ref{fig:2facepolymap}.

Clearly all vertices internal to a given face are at descents of
that face. Therefore, temporarily ignoring  extremal vertices, a
polygon with $j_f-1$ vertices internal to face $f$ (for $f=1,2$)
contributes $u q_{j_1+j_2} z_1^{j_1-1} z_2^{j_2-1}$ to
$\PsmhurGS{2}(z_1,z_2,\vec{q},u)$.  Sum over $j_1, j_2 \geq 1$ to
define
\begin{align*}
  \delta := \sum_{j_1, j_2 \geq 1} u q_{j_1+j_2} z_1^{j_1-1} z_2^{j_2-1}
    = u\frac{Q(z_1)-Q(z_2)}{z_1-z_2}.
\end{align*}
Now the sole extremal vertex incident with polygons $l_{i-1}$ and
$l_i$ is at a descent of face 1 if  $(l_{i-1},l_i)$ is a descent of
$\olist{l_1,\ldots,l_r}$, and at a descent of face 2 otherwise. Thus
Lemma~\ref{lem:circperm} gives
\begin{align*}
    \PsmhurGS{2}(z_1,z_2, \vec{q}, u)
        &=  \log\pr{\frac{x-y}{xe^y-ye^x}}\Bigg|_{x = z_1\delta,\; y = z_2\delta} \\
        &=  \log\pr{\frac{z_1-z_2}{z_1e^{-uQ(z_1)}-z_2e^{-uQ(z_2)}}}
                - u\pr{\frac{z_1 Q(z_1)-z_2 Q(z_2)}{z_1-z_2}}.
\end{align*}
The result follows immediately from Theorem~\ref{thm:cactuspruning}
and identity~\eqref{eq:rcdefn}. \qed

\subsection{Further Results}
\label{ssec:furtherresults}

As mentioned in the introduction, we have been unable to extend the \emph{ad hoc} enumeration applied in our proof
of Theorem~\ref{thm:mainthm1} to give a general formulation of $\PshurGS{m}$ for $m > 3$.  The problem, of course,
is that structure of smooth polymaps on $m$ faces gets far more complicated as $m$ increases; in particular,
polygons can be incident with three or more faces. The case $m=3$ appears within reach (though tedious), and the
more refined enumeration offered by $\PshurGS{3}$ (as compared with Goulden and Jackson's series $\PshurGS{3}{k}$)
may provide some insight into the more general structure.

We have had somewhat more success in applying our techniques to analyze the Hurwitz problem.  Observe that the
polymaps associated with \mbox{2-cycle} factorizations (through Proposition~\ref{prop:mapping}) can be simplified
by collapsing \mbox{2-gons} to single edges, thereby eliminating all polygons and leaving \emph{maps} in the
ordinary sense. In particular, cacti become trees, so \S\ref{ssec:fullcycles} gives a bijection between
\mbox{2-cycle} factorizations of $(1\,2\,\cdots\,n)$ and rooted trees with edges labelled $1,\ldots,n-1$. These
trees, in turn, correspond with trees on $n$ labelled vertices, as can be seen by ``pushing'' edge labels away
from the root onto vertices and assigning label $n$ to the root. Since there are $n^{n-2}$ trees on $n$ labelled
vertices, we have an elegant bijective proof of Hurwitz's formula~\eqref{eq:hurwitz} in the case $m=1$. This
bijection is equivalent to that given by Moszkowski~\cite{moszkowski}.

In fact, the argument used in~\S\ref{ssec:fullcycles} to associate factorizations of full cycles with rooted
cacti applies more generally to give a bijection between cycle factorizations of any $\p \in \Class{\a}$ and
\mbox{$\a$-proper} polymaps with one distinguished descent in each face. We have used this correspondence
(restricted to \mbox{2-cycle} factorizations) together with the pruning bijection established in
Theorem~\ref{thm:cactuspruning} to give the first bijective proofs of~\eqref{eq:hurwitz} when $m=2$ and
$m=3$. (Details can be found in~\cite{irvingphd}.) Thus the simplicity of Hurwitz's formula in these cases is
explained by a corresponding simplicity in the structure of certain smooth maps with two and three faces. We
have as yet been unable to extend our methods to $m \geq 4$, but we hope that the shift in viewpoint afforded
by pruning can be further exploited in this direction.

We should mention that the techniques employed here apply equally well to transitive factorizations
$F=(\s_r,\ldots,\s_1)$ not constrained by the minimality condition.  If $\s_i \in \Class{\b_i}$ and
$\s_r\cdots\s_1 \in \Class{\a}$, then a parity argument applied to~\eqref{eq:genusdefn} shows $nr - \sum_{i=1}^r
\len{\b_i} = n + \len{\a} - 2 + 2g$ for a unique integer $g \geq 0$ called the \emph{genus} of $F$. Thus minimal
transitive factorizations are of genus 0. Unsurprisingly, genus $g$ factorizations correspond with genus $g$
polymaps, and many of the results described here (in particular, Theorem~\ref{thm:cactuspruning} and
Theorem~\ref{thm:inequivpruning} of the next section) have obvious higher genus analogues.  Our focus on minimal
transitive factorizations reflects our belief that understanding the combinatorics in genus 0 is key to an
understanding in all genera.

The graphical approach can also be used to study what we call \emph{$\b$-factorizations}, which are minimal
transitive factorizations $(\rho,\t_r,\ldots,\t_1)$ such that $\rho \in \Class{\b}$ and the factors
$\t_1,\ldots,\t_r$ are transpositions. Thus $[1^n]$-factorizations are synonymous with \mbox{2-cycle}
factorizations.  Determining the number of $\b$-factorizations of a given $\p \in \Class{\a}$ is known as the
\emph{double Hurwitz problem}. It has been the object of recent attention because of known and conjectural links
with geometry, in particular intersection theory; further information can be found in~\cite{gj-vakil}.  Though we
have recovered a number of the results of~\cite{gj-vakil} through combinatorial methods, we have not made further
headway.

\section{Inequivalent Cycle Factorizations}
\label{sec:inequiv}

We now turn to the enumeration of equivalence classes of cycle factorizations under the commutation relation
introduced in \S\ref{ssec:equivalence}.  The graphical methods developed in the previous section will be
reconsidered in this new context, and an appropriate pruning mechanism will be developed. The section ends
with a proof of Theorem~\ref{thm:mainthm2}.

\subsection{Marked and $\a$-Marked Polymaps}
\label{ssec:markedpolymaps}

\newcommand{\imap}{\widetilde{\mathscr{M}}}
Our starting point is the observation that two cycle factorizations are equivalent if and only if (1) they
have precisely the same factors, and (2) the factors moving any given symbol appear in the same order in both
factorizations. Interpreting this through the lens of our graphical correspondence $F \mapsto \amap_F$
(Proposition~\ref{prop:mapping}), we see that permissible commutations of a cycle factorization $F$
correspond with relabellings of the polygons of $\amap_F$ that preserve the relative order of the polygons
incident with any given vertex.

In particular, \emph{two factorizations are equivalent if and only if their corresponding polymaps have the same
descent structure.} Thus the equivalence class $[F]$ containing $F$ is naturally represented by the polymap
$\imap_F$ that results from stripping the polygon labels of $\amap_F$ and instead recording only the location of
descents.

\begin{figure}[t]
\begin{center}
\includegraphics[width=\textwidth]{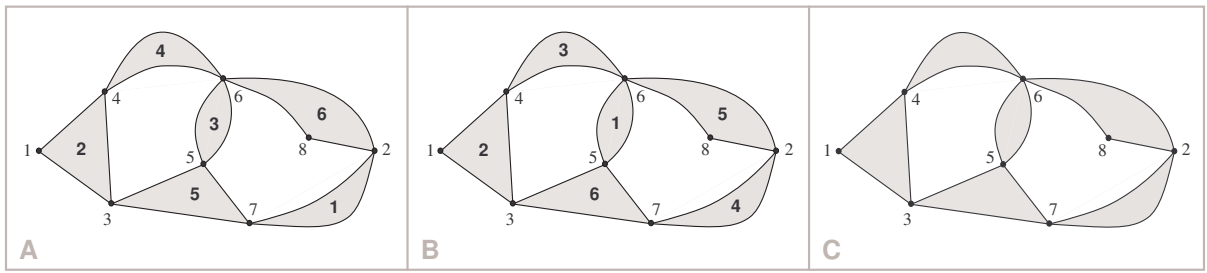}
\caption{Polymaps corresponding to equivalent cycle
factorizations..} \label{fig:inequivalentmaps}
\end{center}
\end{figure}
\begin{exmp} Consider the following equivalent cycle factorizations of $(1\,2\,3)(4\,5)(6\,7\,8)$:
\label{exmp:inequiv}
\begin{align*}
    F &= ((2\,8\,6), (3\,5\,7), (4\,6), (5\,6), (1\,4\,3), (2\,7)) \\
    G &= ((3\,5\,7), (2\,8\,6), (2\,7), (4\,6), (1\,4\,3), (5\,6)).
\end{align*}
Factors $(5\,6)$, $(4\,6)$, and $(2\,8\,6)$ move symbol 6, and they
appear in this same right-to-left order in both factorizations. The
corresponding polymaps $\amap_F$ and $\amap_G$ are shown in
Figures~\ref{fig:inequivalentmaps}A and~\ref{fig:inequivalentmaps}B,
respectively, where vertex labels have been placed in descent
corners. Note that the descent structure of these polymaps is
identical.  The class $[F]$ (equivalently, $[G]$) is represented by
the polymap $\imap_F$ (equivalently, $\imap_G$) shown in
Figure~\ref{fig:inequivalentmaps}C, where again the location of
descents is recorded by the placement of vertex labels. \qed
\end{exmp}

Define a \bold{marked} polymap to be a loopless polymap $\amap$ in
which certain corners have been distinguished subject to the
requirement that every vertex is at exactly one distinguished corner.  A marked polymap
$\amap$ is \bold{properly marked} if its distinguished corners coincide with the descent corners under some
polygon labelling of $\amap$. That is, properly marked
polymaps are exactly those objects that result from marking descents and
deleting polygon labels from proper polymaps. (Note that not every marked polymap is properly marked.)
We refer to the distinguished corners of a properly marked polymap
as \bold{descent corners}, and interpret terms such as \bold{descent
cycle} in the obvious way.

\begin{figure}[t]
    \begin{center}
    \includegraphics[width=\textwidth]{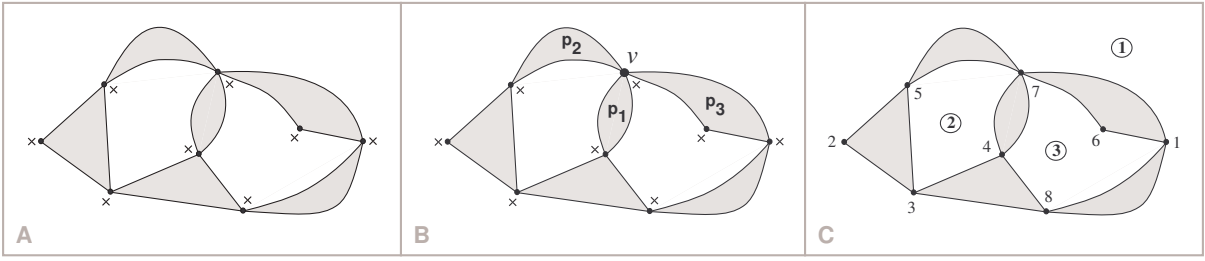}
    \caption{(A,B) A properly marked polymap.  (B) The rotator of $v$ is $(p_1,p_2,p_3)$.  (C) A \mbox{$(3,2,3)$-marked} polymap.}
    \label{fig:markedpolymaps}
    \end{center}
\end{figure}

See Figure~\ref{fig:markedpolymaps}A for an example of a properly marked
polymap.  (Two polygon labellings which validate the indicated descent corners
are shown in Figure~\ref{fig:inequivalentmaps}.)
It will be convenient to redefine the \bold{rotator} of a vertex $v$
in a properly marked polymap to be the tuple $(P_1,\ldots,P_m)$ of polygons
incident with $v$, listed in order as they are encountered along a
clockwise tour about $v$ beginning in the unique descent corner
containing $v$. In Figure~\ref{fig:markedpolymaps}B we have
indicated polygons $p_1,p_2$, and $p_3$ such that the rotator of
vertex $v$ is $(p_1,p_2,p_3)$.

Properly marked polymaps will play the same role in the  enumeration of
inequivalent factorizations to follow as proper polymaps did in our
approach to counting ordered factorizations. Our observations thus
far are summarized by the following analogue of
Proposition~\ref{prop:mapping}:
\begin{prop}
\label{prop:imapping} The mapping $[F] \mapsto \imap_F$ is a bijection between equivalence classes of cycle
factorizations and vertex-labelled properly marked polymaps. If $F$ is a factorization of $\p$ of cycle index $\vec{i}$,
then $\imap_F$ is of poly-index $\vec{i}$ with descent cycles equal to the cycles of $\p$. \qed
\end{prop}

Let $\a=(\a_1,\ldots,\a_m) \cmp n$. An \bold{$\a$-marked} polymap is a vertex-labelled properly marked
polymap whose faces are labelled $1,\ldots,m$ so that face $j$ has descent set $\Dset{\a}{j}$, where
$\Dset{\a}{j}$ is defined as in~\eqref{eq:canonicalsets}. For example, under the convention that the
locations of vertex labels indicate descent corners, the polymap of Figure~\ref{fig:markedpolymaps}C is seen
to be $(3,2,3)$-marked. Of course, the face labels of $\a$-marked polymaps are superfluous, since they are
determined by the descent sets.  They are included in the definition only as a matter of
convenience.\footnote{Note the contrast between this definition of $\a$-marked polymaps and that of
$\a$-proper polymaps in~\S\ref{ssec:properpolymaps}.  In the latter case, \emph{only} face labels were used.
Here we do not have that luxury, because face-labelled marked polymaps can admit nontrivial automorphisms,
even when two or more faces are present.}

Let $\IPshur{\a}{i}$ denote the number of $\a$-marked polymaps with poly-index $\vec{i}$.  Under
Proposition~\ref{prop:imapping}, $\a$-marked polymaps correspond with factorizations of permutations whose
orbits are $\Dset{\a}{1},\ldots,\Dset{\a}{m}$. There are $\prod_j (\a_j-1)!$ such permutations, so
$\IPshur{\a}{i} = \Phur{\a}{i} \cdot \prod_j (\a_j-1)!$. Thus~\eqref{eq:inequivgs} becomes
\begin{equation}
\label{eq:amarkedpolymaps}
   \IPshurGS{m}(\vec{x},\vec{q}, u) =
    \sum_{n \geq 1}    \, \sum_{\vec{i} \,\geq\, \vec{0}}
    \sum_{\substack{\a \cmp n \\ \len{\a} =  m}}
       \IPhur{\a}{\vec{i}} \, \vec{q}^{\vec{i}}\,
        \frac{\vec{x}^{\pmb{\a}}}{\pmb{\a}!} u^{r(\vec{i})}.
\end{equation}
That is, $\IPshurGS{m}(\vec{x},\vec{q},u)$ is the generating series for $\a$-marked polymaps, where $u$ marks
polygons, $\vec{q}$ records poly-index, and $x_j$ is an exponential marker for labelled vertices at descents
of face $j$.

\subsection{Marked Cacti and Inequivalent Factorizations of Full Cycles}
\label{ssec:icacti}

Observe that every marked cactus is  properly marked. The argument used in~\S\ref{ssec:fullcycles} to show
$\Phur{[n]}{i}$ is the number of vertex-rooted proper cacti of poly-index $\vec{i}$ therefore applies
\emph{mutatis mutandis} to identify $\IPhur{[n]}{i}$ as the number of vertex-rooted marked cacti of
poly-index $\vec{i}$. Thus the series $\irt=\irt(x,\vec{p},u)$ defined by
\begin{equation}
\label{eq:ircdefn2}
    \irt := \sum_{n \geq 1} \sum_{\vec{i \geq 0}}
                \IPhur{[n]}{i} x^{n-1} \vec{q}^{\vec{i}}
                u^{r(\vec{i})} =    \frac{d}{dx} \IPshurGS{1}(x,\vec{q},u)
\end{equation}
counts such cacti with respect to the number of \emph{non-root} vertices (marked by $x$), polygons (marked by
$u$), and poly-index (marked by $\vec{q}$). As  in~\S\ref{ssec:fullcycles}, we now consolidate this
definition of $\irt$ with the one given in the introduction by describing a decomposition of marked cacti
that proves their generating series satisfies~\eqref{eq:ircdefn}.

\begin{figure}[t]
    \begin{center}
    \includegraphics[width=\textwidth]{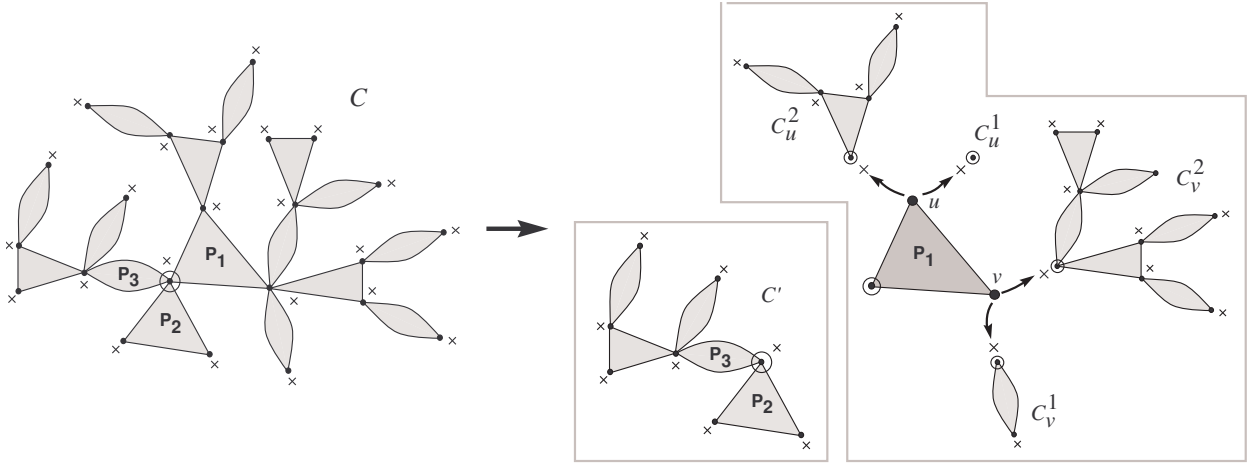}
    \caption{Decomposition of a descent-marked cactus.} \label{fig:inequivalentcactus}
    \end{center}
\end{figure}

Let $C$ be a rooted marked cactus, and suppose its root vertex has
rotator $(p_1,\ldots,p_m)$. Detach $p_1$ from $C$ to obtain a
rooted marked cactus $C'$ whose root has rotator $(p_2,\ldots,p_m)$,
as shown in Figure~\ref{fig:inequivalentcactus}. (If $m=1$ then $C'$
consists of a single vertex.)
Now focus on polygon $p_1$.  Suppose $p_1$ is a $k$-gon, and let $v$
be one of its $k-1$ non-root vertices. Let
$(p_{v}^1,\ldots,p_{v}^r,p_1,p_{v}^{r+1},\ldots,p_{v}^{s})$ be the
rotator of $v$, where the degenerate cases $r=0$ and $s=r$ are
possible. Detach the polygons $p_v^i$ from $v$ to form two
cacti, $C^1_v$ and $C^2_v$, whose roots have rotators
$(p_{v}^1,\ldots,p_{v}^r)$ and $(p_{v}^{r+1},\ldots,p_{v}^{s})$,
respectively. See Figure~\ref{fig:inequivalentcactus}, far right.

Thus $C$ decomposes into $C'$ together with $k$-gon $p_1$
and an ordered list of triples $(v,C^1_v,C^2_v)$, one for each of
the $k-1$ non-root vertices of $p_1$. Allowing for all possible $k$
gives
$$
    \irt = 1 + \sum_{k \geq 2} \irt \cdot u q_k \cdot (x\irt^2)^{k-1}
         = 1+ u\irt Q(x \irt^2),
$$
where the summand 1 accounts for the case in which $C$ consists of a
single vertex, and $Q$ is defined as in~\eqref{eq:qdefn}. This is in
agreement with~\eqref{eq:ircdefn}, as claimed.

The following result first appeared in~\cite{springer}. The bijective proof given there relied on the
commutation of cycle factorizations into canonical forms, together with  a correspondence between these
canonical forms and a class of trees for which enumerative formulae are known. In hindsight, our use of
marked cacti can be viewed as a high level graphical interpretation of this bijection that bypasses the
element-wise analysis necessary to establish the canonical forms.

\begin{thm}[Springer~\cite{springer}]
\label{thm:ispringer} Let $\vec{i}=(i_2,i_3,\ldots)$ be a sequence
of nonnegative integers, not all  zero, and set $r=r(\vec{i})=i_2 +
i_3 + \cdots$. Then the number of inequivalent cycle factorizations
of $(1\,2\,\cdots\,n)$ of cycle index $\vec{i}$ is
$$
    \IPhur{[n]}{i} = \frac{(2n+r-2)!}{(2n-1)! \, \prod_{k \geq 2} i_k!}
$$
in the case that $n + r - 1 = \sum_{k \geq 2} ki_k$, and zero otherwise.
\end{thm}

\begin{proof}Set $v =
\irt-1$.  Then~\eqref{eq:ircdefn2} and~\eqref{eq:ircdefn} give
$\IPhur{[n]}{i}=[x^{n-1} u^r \vec{q}^\vec{i}]\, v$ with $v = u (1+v)
Q(x(1+v)^2).$ By Lagrange inversion~\cite{goulden-jackson} we have
\begin{align*}
    [x^{n-1} u^r \vec{q}^\vec{i}]\,v
    &= [x^{n-1} \vec{q}^\vec{i}]\,\,\frac{1}{r}\, [\l^{r-1}]\,
                (1+\l)^r Q(x(1+\l)^2)^r \\
    &= \frac{1}{r}\, [\l^{r-1}]\,(1+\l)^{r+2n-2}
            \binom{r}{i_2,\, i_3,\, i_4,\,\ldots},
\end{align*}
and the result follows upon simplification.
\end{proof}

Specializing the previous theorem to count inequivalent
$k$-cycle factorizations entails setting $i_k=r$ and $i_j=0$ for $j
\neq k$.  Doing so, we find that the number of inequivalent
$k$-cycle factorizations of $(1\,2\,\cdots\,n)$ is
$$
    \frac{1}{2n-1}\binom{2n+r-2}{r}
$$
when $n=1+r(k-1)$ for some integer $r$.  This formula first appeared
in~\cite{gj-macdonald}.  Setting $k=2$ yields the
result~\eqref{eq:catalan} of Longyear mentioned in the introduction.

\subsection{Pruning Cacti}
\label{ssec:Ipruningcacti}

Let $\a=(\a_1,\ldots,\a_m)$ be a composition.  Define the \bold{face
degree sequence} of an $\a$-marked polymap to be $(d_1,\ldots,d_m)$,
where $d_j$ is the degree of face $j$, for $1 \leq j \leq m$.
Alternatively,  $d_j$ is the number of corners in face $j$. Let
$\IPsmhur{\a}{i}{d}$ denote the number of smooth $\a$-marked
polymaps with poly-index $\vec{i}$ and face degree sequence
$\vec{d}$. Then, for  $m \geq 1$, define
\begin{equation*}
    \IPsmhurGS{m}(\vec{x},\vec{t},\vec{q}, u) :=
    \sum_{n \geq 1}    \, \sum_{\vec{i},\vec{d} \,\geq\, \vec{0}}
    \sum_{\substack{\a \cmp n \\ \len{\a} =  m}}
       \IPsmhur{\a}{i}{d} \,
        \frac{\vec{x}^{\pmb{\a}}}{\pmb{\a}!} \,
                            \vec{q}^{\vec{i}}\, \vec{t}^{\vec{d}}
                            u^{r(\vec{i})}.
\end{equation*}
With these definitions we have the following marked polymap analogue of
Theorem~\ref{thm:cactuspruning}.

\begin{thm}
\label{thm:inequivpruning} Let $m \geq 1$ and set $\irt_i =
\irt(x_i,\vec{q},u)$ for $1 \leq i \leq m$, where $\irt$ is given
by~\eqref{eq:ircdefn}. Then
\begin{equation*}
    \IPshurGS{m}(\vec{x}, \vec{q},  u) =
    \IPsmhurGS{m}(\vec{x}\circ \vec{\irt},\; \vec{\irt},\; \vec{q},\; u),
\end{equation*}
where $\vec{x}=(x_1,\ldots,x_m)$,  $\vec{\irt} = (\irt_1,\ldots,\irt_m)$ and
$\vec{x}\circ\vec{\irt} = (x_1\irt_1,\ldots,x_m\irt_m)$.
\end{thm}

\begin{proof}
Let $\amap$ be a properly marked polymap with at least two faces and let $v$ be a vertex of $\acore$ incident with face $f$.  Let $(p_1,\ldots,p_k)$ be the
rotator of $v$.  We wish to prune from $\amap$ all branches in $f$ attached at $v$.

Suppose first that $v$ is at a descent of $f$.  Then there exist $i,j$ with $1 \leq i \leq j \leq k$ such
that $p_1,\ldots,p_{i-1}, p_{j+1},\ldots,p_k$ are the root polygons of the branches attached to $v$ in face
$f$.  (See Figure~\ref{fig:inequivpruning}A. For clarity, most descent corners have not been indicated.)
Detach these branch  from $\amap$ to form a new polymap $\amap'$ and two rooted marked cacti, $C_1$ and
$C_2$, whose roots have rotators $(p_1,\ldots,p_{i-1})$ and $(p_{j+1},\ldots,p_k)$, respectively. Note that
cases $i=1$ and $j=k$ yield trivial cacti $C_1$ and $C_2$, respectively. The descent corner of $f$ at $v$ is
inherited by $\amap'$, and it is easy to check that this makes $\amap'$ properly marked.
Figure~\ref{fig:inequivpruning}A illustrates this pruning mechanism.
\begin{figure}[t]
\begin{center}
    \psfrag{frag:m}{$\amap$}
    \psfrag{frag:mp}{$\amap'$}
    \includegraphics[width=.85\textwidth]{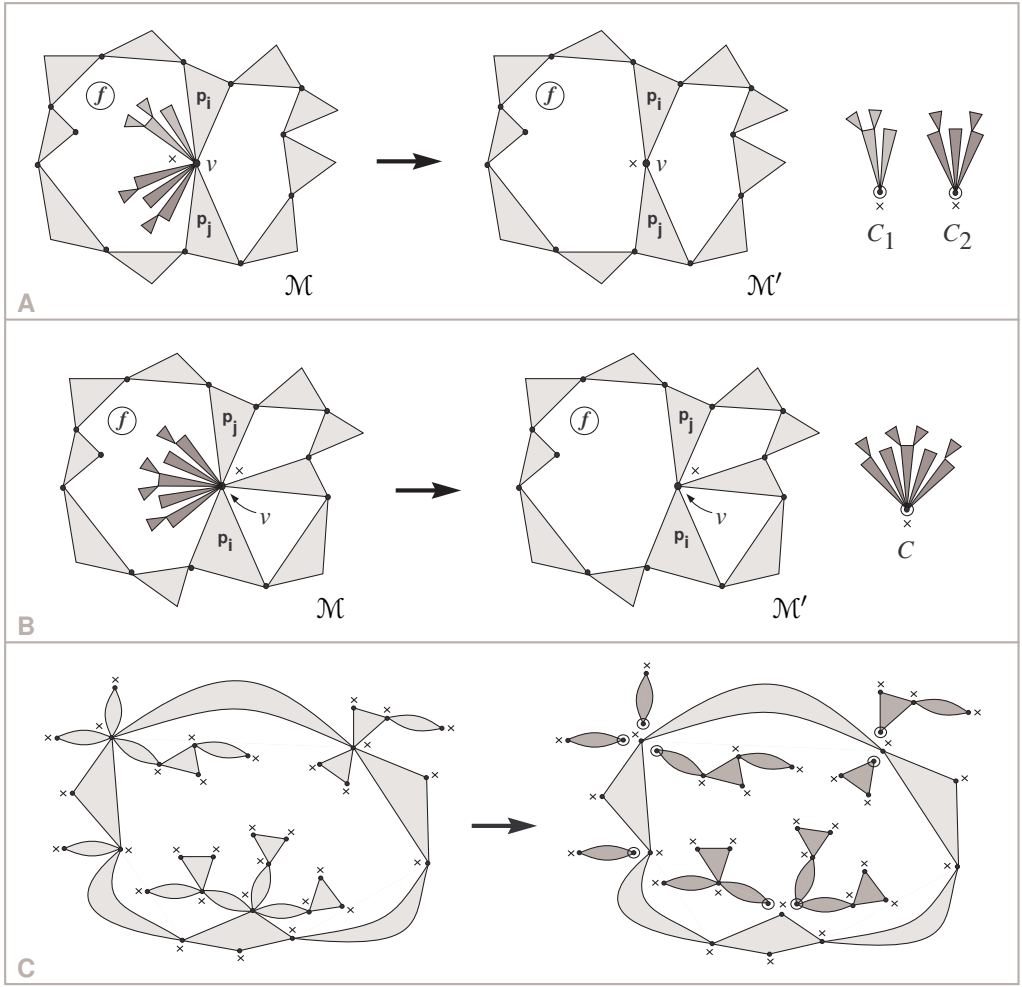}
    \caption{Pruning cacti from a marked polymap.} \label{fig:inequivpruning}
\end{center}
\end{figure}

Now suppose $v$ is not at a descent of $f$.  Then there exist $i,j$ with $1 \leq i < j \leq k$ such that
$p_{i+1},\ldots,p_{j-1}$ are the root polygons of the branches attached to $v$ in $f$.  (See
Figure~\ref{fig:inequivpruning}B.) Detach these branches from $\amap$ to form a  properly marked polymap
$\amap'$ and a rooted marked cactus $C$ whose root has rotator $(p_{i+1},\ldots,p_{j-1})$. When $i=j$ there
are no branches attached to $v$ in $f$, and $C$ is trivial.

Repeating this pruning process for each incident face-vertex pair $(v,f)$ of $\amap$ (that is, for each
corner of $\amap$) results in a smooth properly marked polymap $\amap[S]$ and a collection of rooted marked
cacti, each of which is naturally associated with a particular face of $\amap[S]$. In fact, if a face of
$\amap[S]$ is of degree $D$ and has $d$ descents, then it has $D-d$ non-descent corners and thus exactly
$2d+(D-d)=D+d$ associated pruned cacti. Each non-root vertex of these cacti contributes one descent to the
corresponding face of $\amap$.  The pruning of all cacti from a properly marked polymap is illustrated in
Figure~\ref{fig:inequivpruning}C.

To complete the proof, observe that if $\amap$ is $\a$-marked, where $m=\len{\a} \geq 2$, then its corners
are distinguishable and thus the pruning process described above is reversible. Recall that
$\irt(x,\vec{q},u)$ counts rooted marked cacti, with $x$ recording non-root vertices. Since rooting
eliminates nontrivial automorphisms, these vertices can be freely labelled. From~\eqref{eq:amarkedpolymaps}
it follows that $\IPshurGS{m}(\vec{x},\vec{q},u)$ is obtained from $\IPsmhurGS{m}(\vec{z},\vec{t},\vec{q},u)$
through the substitutions $z_j \mapsto x_j\irt_j$ and $t_j \mapsto \irt_j$, for $1 \leq j \leq m$.
\end{proof}

\subsection{Inequivalent Factorizations of Permutations with Two Cycles}
\label{ssec:Itwofacepolymaps}

We now use the pruning of marked polymaps to prove Theorem~\ref{thm:mainthm2}. As mentioned in the
introduction, this theorem generalizes the main result of~\cite{glj-inequivalent}, and it is worth noting
that the intricate inclusion-exclusion argument used there is avoided entirely by our method. In particular,
the following lemma (which is a routine exercise in exponential generating series) is seen to be a more
natural source of the logarithm appearing in~\eqref{eq:igjresult}.

\begin{lem}
\label{lem:beadlemma} For $m,n \geq 1$, let $c_{n,m}$ be the number
of distinct necklaces made of $n$ labelled white beads and $m$
(independently) labelled black beads, counted up to rotational
symmetry. Then
\begin{equation}
\label{eq:beadeqn}
    \sum_{n,m \geq 1} c_{n,m}  \frac{x^n}{n!}\frac{y^m}{m!} z^{n+m}
    = \log\bigg(1+\frac{xy z^2}{1-z(x+y)}\bigg).
\end{equation}
\qed
\end{lem}

\vspace*{1.5mm}

\noindent\textbf{Proof of Theorem~\ref{thm:mainthm2}:} We begin by determining
$\IPsmhurGS{2}(z_1,z_2,t_1,t_2,\vec{q},u)$.
Let $\amap[S]$ be a smooth \mbox{$\a$-marked} polymap, where $\len{\a}=2$.
Then $\amap[S]$ is a closed chain of polygons, each incident with exactly
two others. As in the proof of Theorem~\ref{thm:mainthm1}, we say a vertex incident with two polygons
is \emph{extremal}.  Since $\amap[S]$ is properly marked, at least one extremal
vertex is at a descent of each face.

\begin{figure}[t]
    \begin{center}
    \includegraphics[width=.55\textwidth]{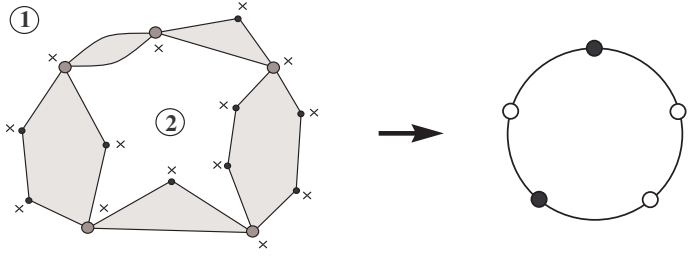}
    \caption{A smooth marked polymap and its associated necklace.}
    \label{fig:maptonecklace}
    \end{center}
\end{figure}

Let $L=\olist{v_1,\ldots,v_r}$ be the cyclic list of extremal
vertices encountered along the boundary walk of face 1 of
$\amap[S]$.  Regard those $v_i$ that are at descents of face 1 as
white beads and those at descents of face 2 as black beads, so $L$
corresponds with a necklace of the type counted by
Lemma~\ref{lem:beadlemma}. (See Figure~\ref{fig:maptonecklace} for
an illustration, where extremal vertices are indicated in grey, but
vertex and bead labels are not shown.)

Let $M$ be the monomial of
$\IPsmhurGS{2}(z_1,z_2,t_1,t_2,\vec{q},u)$ corresponding to
$\amap[S]$. An extremal vertex of $\amap[S]$ contributes the factor
$z_1 t_1 t_2$ to $M$ if it is at a descent of face 1, and
contributes $z_2 t_1 t_2$ otherwise. A polygon of $\amap[S]$ with
$j_f-1$ vertices incident only with face $f$, for $f=1,2$, further
contributes the factor  $u q_{j_1+j_2} (t_1 z_1)^{j_1-1}
(t_2z_2)^{j_2-1}$. Thus $\IPsmhurGS{2}(z_1,z_2,t_1,t_2,\vec{q},u)$
is obtained from~\eqref{eq:beadeqn} by performing the substitutions
$x \mapsto z_1 t_1 t_2$, $y \mapsto z_2 t_1 t_2$, and
\begin{align*}
    z \mapsto \sum_{j_1, j_2 \geq 1} u q_{j_1+j_2} (t_1z_1)^{j_1-1} (t_2z_2)^{j_2-1}
    = u\frac{Q(t_1z_1)-Q(t_2z_2)}{t_1z_1-t_2z_2}.
\end{align*}
But Theorem~\ref{thm:inequivpruning} implies
$\IPshurGS{2}(x_1,x_2,\vec{q},u)$ is got from $\IPsmhurGS{2}(z_1,
z_2, t_1, t_2, \vec{q}, u)$ through the substitutions $z_i \mapsto
x_i \irt_i$ and $t_i \mapsto \irt_i$.  Performing this chain of
substitutions leaves
\begin{align*}
    \IPshurGS{2}(x_1,x_2,\vec{q},u)
    = \log\pr{1 + \frac{\delta^2 x_1 x_2 \irt_1 \irt_2}{1-\delta(x_1 \irt_1 + x_2
    \irt_2)}},
\end{align*}
where
$$
    \delta
    :=
    u\irt_1\irt_2\frac{Q(x_1\irt_1^2)-Q(x_2\irt_2^2)}{x_1\irt_1^2-x_2\irt_2^2}.
$$
Theorem~\ref{thm:mainthm2} now follows upon simplification,
using~\eqref{eq:ircdefn} to write $uQ(x_i\irt_i^2)=1-\irt_i^{-1}$.
\qed

\subsection{Further Results} \label{ssec:ifurtherresults}

\newcommand{\note}[1]{\vspace*{2mm}\noindent\textbf{#1}}

We have been unable to generalize the proof of Theorem~\ref{thm:mainthm2} to obtain expressions for $\IPshurGS{m}$
for any $m > 2$.  However, Theorem~\ref{thm:inequivpruning} has been used~\cite{irvingphd} to give a generating
series formulation for the number of inequivalent \mbox{2-cycle} factorizations of permutations composed of three
cycles. The expression for $\IPkshurGS{3}{2}$ thus obtained is far more complicated than the very
succinct~\eqref{eq:igjresult} and we have been unable to simplify it in a meaningful way.   Finding a natural
expression for $\IPkshurGS{3}{2}$ remains an intriguing open problem, as it would grant far greater insight into
the structure of inequivalent factorizations in general. We note that an approach from the point of view
of~\cite{glj-inequivalent} involves significant technical hurdles that do not surface in the simpler known cases.

Inspired by the definition of $\b$-factorizations (see \S\ref{ssec:furtherresults}) we can extend
Theorem~\ref{thm:ispringer} in the following way. Let $\b = [1^{j_1} 2^{j_2} \cdots] \ptn n$ and consider the set
of minimal transitive factorizations $(\rho,\s_r,\ldots,\s_1)$ of the full cycle $(1\,2\,\cdots\,n)$ such that
$\rho \in \Class{\b}$ and the $\s_i$ are all cycles of length at least two. Say two such factorizations are
equivalent if one can be obtained from the other by commutation of adjacent disjoint \emph{cycle factors}
$\s_1,\ldots,\s_m$. That is, the position of $\rho$ is fixed. Then the number of inequivalent factorizations of
this type with exactly $i_k$ cycle factors being \mbox{$k$-cycles} for $k \geq 2$ is
\begin{equation*}
    \frac{n(r-1)!\,(\len{\b}-1)!}{
        \prod_{k \geq 2} i_k! \cdot \prod_{k \geq 1} j_k!}
    \binom{r+\len{\b}+n-2}{r-1},
\end{equation*}
where $r=i_2+i_3+\cdots$, subject to the necessary condition $\sum
ki_k = r+\len{\b}-1$. This result is readily proved through a
refinement of the decomposition of marked cacti given in
\S\ref{ssec:icacti}. Note that Theorem~\ref{thm:ispringer} is
recovered by setting $\b=[1^n]$.

We conclude by noting that our definition of equivalence of
factorizations may not seem the most natural.  For instance, one
could allow all pairs of adjacent commuting factors to be
interchanged, whether or not they are disjoint.  A general statement
of the corresponding enumerative problem might then be the
following: For a multiset $B$ of partitions of $n$, find the number
of inequivalent minimal transitive factorizations of a given $\p \in
\Sym{n}$ whose factors have cycle types specified by $B$. We have
not yet investigated this problem in any detail.

\section*{Acknowledgements}

The bulk of the research for this article was completed during the author's doctoral studies under the
supervision of David Jackson, whose support and guidance has been greatly appreciated.


\providecommand{\bysame}{\leavevmode\hbox to3em{\hrulefill}\thinspace}
\providecommand{\MR}{\relax\ifhmode\unskip\space\fi MR }
\providecommand{\MRhref}[2]{%
  \href{http://www.ams.org/mathscinet-getitem?mr=#1}{#2}
}
\providecommand{\href}[2]{#2}


\end{document}